\documentclass[11pt]{article}
\usepackage{amsmath,amssymb,amsfonts,amsthm}
\usepackage{latexsym}
\usepackage[english]{babel}
\usepackage{showidx}
\usepackage{graphicx}
\usepackage{graphics}
\usepackage{babel,amsmath,amssymb,amsbsy,amsfonts,latexsym,amsthm,mathrsfs}

\def\<{\langle}
\def\>{\rangle}

\def\card{{\mathrm{card}}}

\def\I{\mbox{\large \bf 1}}
\def\P{{\mathbb P}}
\def\R{{\mathbb R}}
\def\E{{\mathbb E}}
\def\D{{\mathbb D}}
\def\N{{\mathbb N}}
\def\leb{\mathrm{Leb}}
\def\Var{\mathrm{Var}}

\def\SS{\mathcal{S}}
\def\mD{\mathcal{D}}
\def\mK{\mathcal{K}}
\def\mA{\mathcal{A}}
\def\mE{\mathcal{E}}
\def\mH{\mathcal{H}}
\def\mM{\mathcal{M}}

\def\coeffbin#1#2{\Big(\!\begin{array}{c}#1\cr #2\end{array}\!\Big)}

\numberwithin{equation}{section}

\newtheorem{theorem}{Theorem}[section]

\newtheorem{definition}[theorem]{Definition}

\newtheorem{lemma}[theorem]{Lemma}

\newtheorem{proposition}[theorem]{Proposition}
\newtheorem{remark}[theorem]{Remark}

\begin{document}

\parindent 0pt
\title{\textbf{Asymptotic development for the CLT\\in total variation distance}}
\author{ \textsc{Vlad Bally}\thanks{%
Universit\'e Paris-Est, LAMA (UMR CNRS, UPEMLV, UPEC), MathRisk INRIA, F-77454
Marne-la-Vall\'{e}e, France. Email: \texttt{bally@univ-mlv.fr} }\smallskip \\
\textsc{Lucia Caramellino}\thanks{%
Dipartimento di Matematica, Universit\`a di Roma - Tor Vergata, Via della
Ricerca Scientifica 1, I-00133 Roma, Italy. Email: \texttt{%
caramell@mat.uniroma2.it}. Corresponding author.}\smallskip\\
}
\maketitle



\begin{abstract}
The aim of this paper is to study the asymptotic expansion in total variation in the
central limit theorem when the law of the basic random variable is locally
lower-bounded by the Lebesgue measure (or equivalently, has an absolutely continuous component): we develop  the error in powers
of $n^{-1/2}$ and give an explicit formula for the approximating measure.
\end{abstract}


\medskip


\section{Introduction}

The aim of this paper is to study the convergence in total variation in the
Central Limit Theorem (CLT) under a certain regularity condition for the random variable at hand.
Given two measures $\mu,\nu$ in $\R^N$, we recall that the distance in total variation is defined as
$$
d_{TV}(\mu,\nu)=\sup\Big\{\Big|\int fd\mu-\int f d\nu\Big|\,:\,\|f\|_\infty\leq 1\Big\}.
$$
Let $F$ be a centred r.v. in $\R^N$ with identity covariance matrix and let $F_k$, $k\in\N$, denote independent copies of $F$. We set
$$
S_n=\frac 1{\sqrt n}\sum_{k=1}^nF_k.
$$
We also define $\mu_n$ the law of $S_n$ and $\Gamma$ the standard Gaussian law in $\R^N$.

The problem of the convergence in total variation for the CLT, that is $d_{TV}(\mu_n,\Gamma)\to 0$ as $n\to\infty$, is very old. Prohorov \cite{[P]} in 1952 proved that, in dimension 1, a necessary and sufficient condition in order to get the result is that there exists $n_0$ such that the law of $\sum_{k=1}^{n_0}F_k$ has an absolutely continuous component (see next Definition \ref{acc}). Then many related problems have been considered in the literature, such as the generalization to the multidimensional case, the study of the speed of convergence, the convergence and the development of the density of $S_n$, if it exists, or the case of a r.v. $F$ whose law has not necessarily an absolutely continuous component, the latter implying the use of a different distance, which is similar to the total variation one but defined on a special class of test functions, typically indicator functions of special sets.

A first class of results has been obtained by Rao \cite{[R]} and then improved by Battacharaya \cite{[B]}:
in \cite{[R]} one proves that the convergence in the CLT holds when the test
function is the indicator function of a convex set $D$. This result is
improved in \cite{[B]} where $D$ is no more a convex set but a set with a boundary
which is small in some sense. An overview on this topic is given in \cite{[BR]}.
But it turns out that one is not generally able to extend the above
mentioned results to a general set $D$ (and so to general measurable and
bounded test functions), because, thanks to the Prohorov's result, one needs to assume a little bit
of regularity on the law of the basic random variable $F$ which comes on in
the CLT. In such a case, Sirazhdinov and Mamatov \cite{[MS]} prove that if $F\in L^3(\Omega)$ then the density of the absolutely continuous component of the law $\mu_n$ converges in $L^1(\R^N)$ to the standard Gaussian density and, therefore, the convergence of the CLT holding in total variation distance, at speed $1/\sqrt n$. This is done in the one-dimensional case, but it works as well in the multidimensional case. The second part of the book \cite{[BR]} gives a complete view on the recent research on this topic, mainly on the the development of the density of $S_n$ around the standard Gaussian density. Results concerning the convergence in the entropy distance (under the same type of hypothesis) has
been recently obtained in \cite{[BCG]}.

This paper contributes in this direction by giving the precise expansion of the CLT in total variation distance. More precisely, we assume that the law of $F$ is locally lower bounded by the Lebesgue measure $\leb_{N}$ on $\R^{N}$ in the following sense: there exists an open set $D_0$ and $\varepsilon_0 >0$ such that for every
Borel set $A$ one has
\begin{equation}
\P(F\in A)\geq \varepsilon_0 \times \leb_{N}(A\cap D_0).  \label{i1}
\end{equation}%
We will show that this is equivalent to the request that the law of $F$ has an absolutely continuous component (and moreover, we can construct such absolutely continuous measure in order that the associated density is a non-negative lower semicontinuous function, see Appendix \ref{app-1}). So it
is clear that our hypotheses overlaps the assumption of the existence of the density but one cannot
reduce one to another (if the law of $F$ gives positive probability to the rational points then it is not absolutely continuous; and doing convolutions does
not help). Let us give a non-trivial example. Consider a functional $F$ on
the Wiener space and assume that $F$ is twice differentiable in Malliavin
sense: $F\in \D^{2,p}$ with $p>N$ where $N$ is the dimension of $F$. Let $%
\sigma _{F}$ be the Malliavin covariance matrix of $F.$ If $\P(\det \sigma
_{F}>0)=1$ then the celebrated criterion of Bouleau and Hirsh ensures that
the law of $F$ is absolutely continuous, so we are in the classical case (in
fact it suffices that $F\in \D^{1,2}).$ But if $\P(\det \sigma _{F}>0)<1$ this
criterion does no more work (and one may easily produce examples when the
law of $F$ is not absolutely continuous). In \cite{[BC]}, we proved
that if $\P(\det \sigma _{F}>0)>0$ then the law of $F$ has the property (\ref%
{i1}). Notice also that in the one-dimensional case ($N=1)$ the fact that $F$
is not constant immediately implies that $\P(\sigma _{F}>0)>0.$ Indeed, in
this case $\sigma _{F}=\left\vert DF\right\vert ^{2}$ and if this is almost
surely null, then $F$ is constant.

Let us introduce our results. We consider a random variable $F\in L^{2}(\R^{N})$
which satisfies (\ref{i1}), such that $\E(F)=0$ and the covariance matrix of $F$ is the identity matrix. We take a sequence $F_{k},k\in \N$ of independent
copies of $F$ and we denote by $\mu _{n}$ the law of $S_n=\frac{1}{n^{1/2}}%
\sum_{k=1}^{n}F_{k}$ and by $\Gamma $ the standard Gaussian law on $\R^{N}.$
Under these hypotheses, we first prove that $\lim_{n\rightarrow \infty }d_{TV}(\mu
_{n},\Gamma )=0$ where $d_{TV}$ is the total variation distance. Then we
give the asymptotic development, which we are able to find according to additional requests on the existence of the moments of $F$. More precisely, we get that, for $r\geq 2$, if $F\in L^{r+1}(\Omega)$  and if the moments of $F$ up to order $r$  agree with the moments of the standard Gaussian law then under (\ref{i1}) one has
$$
d_{TV}(\mu_n, \Gamma)\leq C
(1+\E(\left\vert F\right\vert ^{r+1}))^{[r/3]\vee 1}\times \frac{1}{n^{(r-1)/2}}.
$$
In the general case, we obtain the following asymptotic expansion.
For $r\geq 2$ and $n\geq 1$, we define a measure on $\R^N$ through
\begin{equation}\label{Gammanp-1-intro}
\Gamma_{n,r}(dx)=\gamma(x)\Big(1+\sum_{m=1}^{[r/3]} \frac 1{n^{\frac m2}}\mK_m(x)\Big)dx,
\end{equation}
where $\gamma$ denotes the probability density function of a standard normal random variable in $\R^N$ and $\mK_m(x)$ is a polynomial of order $m$ ($[\cdot]$ standing for the integer part). Note that for $r=2$ one gets $\Gamma_{n,r}(dx)=\gamma(x) dx=\Gamma(dx)$. So, we prove that if $F\in L^{r+1}(\Omega)$ with $r\geq 2$ then there exist polynomials $\mK_m(x)$, $m=1,\ldots,[r/3]$ (no polynomials are needed for $r=2$), such that, setting $\Gamma_{n,r}$ the measure in (\ref{Gammanp-1-intro}) and $\mu_n$ the law of $S_n$, under (\ref{i1}) one has
\begin{equation}\label{dev-intro}
d_{TV}(\mu_n,\Gamma_{n,r})\leq C(1+\E(|F|^{r+1}))^{[r/3]\vee 1+1}\times \frac{1}{n^{([r/3]+1)/2}}
\end{equation}
where  $C>0$ depends on $r$ and $N$. So, in order to improve the development (and the rate of convergence) one needs to pass from the request $F\in L^{3k}$ to $F\in L^{3k+3}$, $k\geq 1$.

The development given in (\ref{dev-intro}) is analogous to the one obtained in Theorem 19.5, page 199 in \cite{[BR]}. But our development is explicit: in \cite{[BR]}, the result is obtained using the Fourier transform and consequently the coefficients in the development involve the inverse of the Fourier transform, whereas here we give an explicit expression for the polynomials $\mK_m(x)$, as a linear combination of the Hermite polynomials (see next formula (\ref{Dbis})).

The main instrument used in this paper is the Malliavin-type finite
dimensional calculus defined in \cite{[BCl]} and \cite{[BC]}. It turns out that for a
random variable which satisfies (\ref{i1}) a very pleasant calculus may
be settled. The idea is that (\ref{i1}) guarantees that the law of $F$
contains some smooth noise. Then, using a splitting procedure (see Proposition \ref{prop-chi} for details), we may
isolate this noise and achieve integration by parts formulae based on it.

In the last years, a number of results concerning the weak convergence of
functionals on the Wiener space using Malliavin calculus and Stein's method
have been obtained by Nurdin, Peccati, Nualart and Poly; see, for example, \cite{[NNPy],[NP],[NPy],[NPy2]}.
In particular, in
\cite{[NNPy], [NP]} the authors consider functionals living in a finite direct
sum of chaoses and prove that under a very weak non-degeneracy condition
(analogous  to the one we consider here) the convergence in distribution of a
sequence of such functionals implies the convergence in total variation. The
results proved in these papers may be seen as variants of the CLT but for
dependent random variables -- so the framework and the arguments are rather
different from the one considered here.

\section{Main results}

Let $X$ be a random variable in $\R^N$ and let $\mu_X$ denote its law. The Lebesgue decomposition of $\mu_X$ says that there exist a measure $\mu(dx)=\mu(x)dx$, that is,  $\mu$ is absolutely continuous w.r.t. the Lebesgue measure, and a further measure $\nu$ which is singular, that is,  concentrated on a set of null Lebesgue measure, such that
\begin{equation}\label{dec-Leb}
\mu_X(dx)=\mu(x)dx+\nu(dx).
\end{equation}

\begin{definition}\label{acc}
$X$ is said to have an $\mathrm{absolutely\  continuous \ component}$ if the absolutely continuous measure $\mu$ in the decomposition (\ref{dec-Leb}) is not null, that is,  $\nu(\R^N)<1$.
\end{definition}

Definition \ref{acc} plays a crucial role when dealing with the convergence of the Central Limit Theorem (CLT) in the total variation distance $d_{TV}$. We recall the definition of $d_{TV}$: for any two measures $\mu$ and $\nu$ in $\R^N$ then
$$
d_{TV}(\mu,\nu)=\sup\Big\{\Big|\int fd\mu-\int fd\nu\Big|\,:\,\|f\|_\infty\leq 1\Big\},
$$

We discuss here the CLT in total variation distance, so we consider a sequence $\{F_k\}_k$  of i.i.d. square integrable random variables, with null mean and covariance matrix $C(F)$. We set $A(F)$ the inverse of $C(F)^{1/2}$ and
$$
S_n=\frac{1}{\sqrt n}\sum_{k=1}^n A(F)F_k.
$$
We recall the following classical result, due to Prohorov \cite{[P]}.

\begin{theorem}\label{prohorov}
$\mathrm{[Prohorov]}$
Let $\mu_n$ denote the law of $S_n$ and $\Gamma$ denote the standard Gaussian law in $\R^N$. The convergence in the CLT takes place w.r.t. the total variation distance, that is $d_{TV}(\mu_n,\Gamma)\to 0$ as $n\to\infty$, if and only if there exists $n_0 \geq 1$
such that the random variable $S_{n_0}$ has an absolutely continuous component.
\end{theorem}

Hereafter, we assume that the common law of the $F_k$'s has an absolutely continuous component, and this is not a big loss in generality.  In fact, due to the Prohorov's theorem, otherwise we can packet the sequence $\{F_k\}_k$ in groups of $n_0$ r.v.'s, so we can deal with
$$
\bar S_{n}=\frac 1{\sqrt n}\sum_{k=1}^n\bar F_k\quad\mbox{where}\quad \bar F_k=\frac1{\sqrt n_0}\sum_{i=kn_0}^{(k+1)n_0}F_i.
$$
Let us introduce an equivalent way to see probability laws having an absolutely continuous component. From now on, $\leb_N$ denotes the lebesgue measure in $\R^N$.

\begin{definition}\label{acc-bis}
A probability law $\mu$ in $\R^N$ is said to be {\rm locally lower bounded by the Lebesgue measure}, in symbols $\mu\succeq \leb_N$, if there exist $\varepsilon_0>0$ and an open set $D_0\subset \R^N$ such that
\begin{equation}
\mu(A)\geq \varepsilon_0 \leb_{N}(A\cap D_0)\qquad \forall A\in
\mathcal{B}(\R^{N}).  \label{abs1}
\end{equation}%
\end{definition}

We have the following.

\begin{proposition}\label{prop-main1}
Let $F$ be a r.v. in $\R^N$ and let $\mu_F$ denote its law. Then the following statements are equivalent:

 \smallskip

$(i)$ $\mu_F\succeq \leb_N$;

 \smallskip

$(ii)$
$F$ has an absolutely continuous component;

\smallskip

$(iii)$
there exist three independent r.v.'s $\chi$  taking values in $\{0,1\}$, with $\P(\chi=1)>0$,  and $V,W$ in
$\R^{N}$, with $V$ absolutely continuous, such that
\begin{equation}
\P(\chi V+(1-\chi )W\in dv)=\mu _{F}(dv).  \label{abs7}
\end{equation}%

\medskip

Moreover, if one of the above conditions holds then the covariance matrix $C(F)$ of $F$ is invertible.

\end{proposition}

The proof of Proposition \ref{prop-main1} is postponed to Appendix \ref{app-1}. As an immediate consequence of Proposition \ref{prop-main1},  if $\mu _{F}\succeq \leb_{N}$ then $\underline{\lambda}(F)>0$, $\underline{\lambda}(F)$ denoting the smallest eigenvalue of $\widehat{C}(F)=C(F)^{-1}$. We denote through  $\overline{\lambda}(F)$ the associated largest eigenvalue.

\smallskip

We are now ready to introduce the main contributions of this paper. We first give a new proof of the convergence in total variation in the CLT.

\begin{theorem}
\label{TV}Suppose that $\mu _{F}\succeq \leb_{N},$ $\E(F)=0$ and $\E(\left\vert
F\right\vert ^{2})<\infty .$ Then%
\begin{equation}
\lim_{n\rightarrow \infty }d_{TV}(\mu_{n},\Gamma )=0  \label{abs22}
\end{equation}%
where $\mu_n$ denotes the law of $S_n$ and $\Gamma$ is the standard Gaussian law in $\R^N$.
\end{theorem}

This is done especially in order to set up the main arguments and results from abstract Malliavin calculus coming from representation (\ref{abs7}), that are used throughout this paper. Let us stress that Nourdin and Poly in \cite{[NPy2]} have dealt with r.v.'s fulfilling properties that imply (\ref{abs7}), to which they apply results from \cite{[BGL]} about a finite dimensional Malliavin type calculus.

Afterward, we deal with the estimate of the error. In fact, by means of additional requests of the existence of the moments of $F$ up to order $\geq 3$, we get the asymptotic expansion in powers of $n^{-1/2}$ of the law of $S_n$ in total variation distance.
We first obtain the following.
\begin{theorem}
\label{Speed-r}
Suppose that $\mu _{F}\succeq \leb_{N}$ and $\E(F)=0$. Let $\mu_n$ denote the law of $S_n$ and $\Gamma$ denote the standard Gaussian law in $\R^N$. Let $r\geq 2$.
If $\E(\left\vert F\right\vert ^{r+1})<\infty$ and all moments up to order $r$  of $A(F)F$ agree with the moments of a standard Gaussian r.v. in $\R^N$ then
\begin{equation}\label{abs35''bis}
d_{TV}(\mu_n, \Gamma)\leq C
(1+\E(\left\vert F\right\vert ^{r+1}))^{[r/3]\vee 1}\times \frac{1}{n^{(r-1)/2}}
\end{equation}
where $C>0$ depends on $r$, $N$, $\underline{\lambda}(F)$ and $\overline{\lambda}(F)$.

\end{theorem}

In the general case, that is the moments do not generally coincide, we get the following expansion.
For $r\geq 2$ and $n\geq 1$, we define a measure on $\R^N$ through
\begin{equation}\label{Gammanp-1}
\Gamma_{n,r}(dx)=\gamma(x)\Big(1+\sum_{m=1}^{[r/3]} \frac 1{n^{\frac m2}}\mK_m(x)\Big)dx,
\end{equation}
where $\gamma$ denotes the probability density function of a standard normal random variable in $\R^N$ and $\mK_m(x)$ is a polynomial of order $m$ -- the symbol $[\cdot]$ stands for the integer part and for $r=2$   the sums in (\ref{Gammanp-1}) nullify, so that $\Gamma_{n,2}(dx)=\gamma(x) dx=\Gamma(dx)$. Then we get the following.

\begin{theorem}\label{main-th-1}
Let $r\geq 2$ and $\E(|F|^{r+1})<\infty$. Then there exist polynomials $\mK_m(x)$, $m=1,\ldots,[r/3]$ (no polynomials are needed for $r=2$), such that, setting $\Gamma_{n,r}$ the measure in (\ref{Gammanp-1}) and $\mu_n$ the law of $S_n$, one has
$$
d_{TV}(\mu_n,\Gamma_{n,r})\leq C(1+\E(|F|^{r+1}))^{[r/3]\vee 1}\times \frac{1}{n^{([r/3]+1)/2}}
$$
where $C>0$ depends on $r$, $N$, $\underline{\lambda}(F)$ and $\overline{\lambda}(F)$.
\end{theorem}

The statement of Theorem \ref{main-th-1} is not properly written, because no information is given about the polynomials $\mK_m$'s. We observe  that in next formula (\ref{Dbis}) we give a closed-form expression for the $\mK_m$'s in terms of a linear combination of Hermite polynomials, whose coefficients can be explicitly written (so not involving inverse Fourier transforms).

\begin{remark}\label{mall}
Let $F\in\D^{2,p}$ with $p>N$, $\D^{k,p}$ denoting the set of the random variables which are derivable in Malliavin sense up to order $k$ in $L^p$ (see Nualart \cite{[N]}).  If $\P(\sigma_F>0)>0$, $\sigma_F$ standing for the Malliavin covariance matrix of $F$ (and note that this request is much weaker than the non-degeneracy of $\sigma_F$) then Theorem 2.16 in \cite{[BC]} gives that $\mu_F\succeq\leb_N$ (and this property may be strict, that is $F$ may not be absolutely continuous). So both Theorem \ref{Speed-r} and Theorem \ref{main-th-1} can be applied.
\end{remark}

The rest of this paper is devoted to the proofs of the above results: Section \ref{sect-TV} allows us to prove Theorem \ref{TV} and the remaining Theorem \ref{Speed-r} and Theorem \ref{main-th-1} are discussed in Section \ref{sect-asympt}.

\section{Convergence in the total variation distance}\label{sect-TV}

The aim of this section is to prove Theorem \ref{TV}, whose proof requires some preparatives which will be useful also in the sequel.

\subsection{Abstract Malliavin calculus based on a splitting method}\label{AMC}

We consider a random variable $F\in \R^{N}$ whose law $\mu
_{F} $ is such that $\mu_F\succeq\leb_N$. As proved in  Proposition \ref{prop-main1}, the covariance matrix $C(F)$ of $F$ is invertible. So, without loss of generality we can assume from now on that $C(F)$ is the identity matrix, otherwise we work with $A(F)F$, $A(F)$ being the inverse of $C(F)^{1/2}$.

We consider the following special splitting for the law of $\mu_F$, giving, as a consequence, representation (\ref{abs7}). We start from the class of localization functions $\psi _{a}:{%
\mathbb{R}}\rightarrow {\mathbb{R}}$, $a>0$, defined as
\begin{equation}
\psi _{a}(x)=1_{|x|\leq a}+\exp \Big(1-\frac{a^{2}}{a^{2}-(|x|-a)^{2}}\Big)%
1_{a<|x|<2a}.  \label{abs2}
\end{equation}%
Then $\psi _{a}\in C_{c}^{\infty }({\mathbb{R}})$ (the subscript ``$c$'' standing for compact support), $0\leq \psi _{a}\leq 1$
and we have the following property: for every $k,p\in \N$ there exists a
universal constant $C_{k,p}$ such that for every $x\in \R_{+}$%
\begin{equation}
\psi _{a}(x)\left\vert (\ln \psi _{a})^{(k)}(x)\right\vert ^{p}\leq \frac{%
C_{k,p}}{a^{pk}}.  \label{abs3}
\end{equation}%
By the very definition, if $\mu_F\succeq\leb_N$ then we may find $v_{0}\in \R^{N},r_0>0$ and $\varepsilon_0 >0$
such that $\P(F\in A)\geq \varepsilon_0 \leb_{N}(A\cap B_{r_0}(v_{0})).$ Then for
every non-negative function $f:\R^{N}\rightarrow \R_{+}$ we have%
\begin{equation}
\E(f(F))\geq \varepsilon_0 \int_{\R^{N}}\psi _{r_0/2}(\left\vert
v-v_{0}\right\vert )f(v)dv.  \label{abs4}
\end{equation}%
We denote
\begin{equation}
m_0=\varepsilon_0 \int_{\R^{N}}\psi _{r_0/2}(\left\vert v-v_{0}\right\vert )dv.
\label{abs5}
\end{equation}%
Of course, $m_0>0$. But, up to choose $\varepsilon_0$ smaller, we also have $m_0<1$.
So, we consider three independent random variables $\chi \in \{0,1\}$ and $V,W\in
\R^{N}$ with laws
\begin{equation}\label{ans6}
\begin{array}{c}
\displaystyle
\P(\chi  =1)=m_0,\qquad\quad \P(\chi =0)=1-m_0,   \smallskip\\
\displaystyle
\P(V \in dv)=\frac{\varepsilon_0}{m_0}\, \psi _{r_0/2}(\left\vert
v-v_{0}\right\vert )dv, \smallskip\\
\displaystyle
\P(W \in dv)=\frac{1}{1-m_0}\big(\mu _{F}(dv)-\varepsilon_0 \psi _{r_0/2}(\left\vert
v-v_{0}\right\vert )dv\big).
\end{array}%
\end{equation}
Then
\begin{equation}
\P(\chi V+(1-\chi )W\in dv)=\mu _{F}(dv).  \label{abs7-OLD}
\end{equation}%
So, we have just proved the following
\begin{proposition}\label{prop-chi}
If $\mu_F\succeq\leb_N$ then representation (\ref{abs7}) holds.
\end{proposition}

From now on we will work with the representation of $\mu _{F}$ in (\ref{abs7-OLD}) so we always
take
$$
F=\chi V+(1-\chi )W,
$$
$\chi$, $V$ and $W$ being independent and whose laws are given in (\ref{ans6}).

We come now to the central limit theorem. We consider a sequence $\chi
_{k},V_{k},W_{k}\in \R^{N},k\in \N$ of independent copies of $\chi ,V,W\in
\R^{N}$ and we take $F_{k}=\chi _{k}V_{k}+(1-\chi _{k})W_{k}.$ Then we look
to
\begin{equation*}
S_{n}=\frac{1}{n^{1/2}}\sum_{k=1}^{n}F_{k}
=\frac{1}{n^{1/2}}\sum_{k=1}^{n}\big(\chi _{k}V_{k}+(1-\chi _{k})W_{k}\big).
\end{equation*}%
In order to prove the CLT in the total variation distance, we
will use the abstract Malliavin calculus settled in \cite{[BCl]}
and \cite{[BC]} associated to the basic noise
\begin{equation}\label{V}
V=(V_{1},\ldots,V_{n})=((V_{1}^{1},...,V_{1}^{N}),\ldots,(V_{n}^{1},...,V_{n}^{N}))%
\in \R^{N\times n}
\end{equation}
(this will be done for each fixed $n$). To begin, we
recall the notation and some results from \cite{[BC]}. We work with
functionals $X=f(V)$ with $f\in C_{b}^{\infty }(\R^{N\times n};\R)$, the subscript ``$b$'' standing for bounded derivatives of any order. Then we set
$$
\SS=\{f(V)\,:\,f\in C_{b}^{\infty }(\R^{N\times n};\R)\}
$$
and for a functional $X\in \SS$ we define the Malliavin derivatives
\begin{equation}
D_{(k,i)}X=\frac{\partial X}{\partial V_{k}^{i}}=\frac{\partial f}{\partial
v_{k}^{i}}(V),\qquad k=1,...,n,i=1,...,N.  \label{abs9}
\end{equation}%
The Malliavin covariance matrix for a multidimensional functional $%
X=(X^{1},...,X^{d})\in\SS^d$ is defined as
\begin{equation}
\sigma _{X}^{i,j}=\left\langle DX^{i},DX^{j}\right\rangle
=\sum_{k=1}^{n}\sum_{r=1}^{N}D_{(k,r)}X^{i}\times D_{(k,r)}X^{j},\quad i,j=1,\ldots,d.
\label{abs13}
\end{equation}%
We will denote by $\lambda _{X}$ the lower eigenvalue of $\sigma _{X}$, that
is,
\begin{equation}
\lambda _{X}=\inf_{\left\vert \xi \right\vert =1}\left\langle \sigma _{X}\xi
,\xi \right\rangle =\inf_{\left\vert \xi \right\vert
=1}\sum_{k=1}^{n}\sum_{i=1}^{N}\left\langle D_{(k,i)}X,\xi \right\rangle
^{2}.  \label{abs13a}
\end{equation}%
Moreover we define the higher order derivatives just by iterating $D$. We consider a
multiindex $\alpha =(\alpha _{1},...,\alpha _{m})$ with $\alpha
_{j}=(k_{j},i_{j}),k_{j}\in \{1,...,n\},i_{j}\in \{1,...,N\}$ and we set $|\alpha|=m$. Then, we define
\begin{equation}
D_{\alpha }X=\frac{\partial ^{m}X}{\partial V_{k_{m}}^{i_{m}}....\partial
V_{k_{1}}^{i_{1}}}=\partial _{\alpha }f(V)  \label{abs10}
\end{equation}%
with%
\begin{equation*}
\partial _{\alpha }f(v)=\frac{\partial ^{m}f}{\partial
v_{k_{m}}^{i_{m}}....\partial v_{k_{1}}^{i_{1}}}(v).
\end{equation*}%
We will work with the norms%
\begin{eqnarray}
\left\vert X\right\vert _{1,m}^{2} &=&\sum_{1\leq \left\vert \alpha
\right\vert \leq m}\left\vert D_{\alpha }X\right\vert ^{2},\qquad \left\vert
X\right\vert _{m}^{2}=\left\vert X\right\vert ^{2}+\left\vert X\right\vert
_{1,m}^{2}  \label{abs11} \\
\left\Vert X\right\Vert _{1,m,p} &=&\| \,| X|
_{1,m}\,\|_{p}=(\E(\left\vert X\right\vert _{1,m}^{p}))^{1/p},\qquad
\left\Vert X\right\Vert _{m,p}=\left\Vert X\right\Vert _{p}+\left\Vert
X\right\Vert _{1,m,p}  \label{abs11'}
\end{eqnarray}%
We define now the Ornstein--Uhlenbeck operator by%
\begin{equation}
-LX=\sum_{k=1}^{n}\sum_{i=1}^{N}D_{(k,i)}D_{(k,i)}X+\sum_{k=1}^{n}%
\sum_{i=1}^{N}D_{(k,i)}X\partial _{i}\ln \psi _{r_0/2}(\left\vert
V_{k}-v_{0}\right\vert ).  \label{abs12}
\end{equation}%
These are the operators introduced in \cite{[BCl]} and \cite{[BC]} in connection to
the random variable $V$ in (\ref{V}) and taking the weights $\pi _{k}=1$. We will use the
results from \cite{[BC]} in this framework. In particular, as a straightforward consequence of Theorem 3.1 in \cite{[BC]} (take $\mathbf{\Theta}=1$ therein) and  Theorem 3.4 in \cite{[BC]} (see (3.28) therein), we can state integration by parts formulas and estimates for the weights. For later use, we resume in the following statement such facts:

\begin{proposition}
\label{IBP} $X\in \mathcal{S}^{d}$ be such that
$$
\|(\det \sigma_{X})^{-1}\|_p<\infty\quad\mbox{for every}\quad p\geq 1.
$$
Set  $\gamma _{X}$ the inverse of $\sigma _{X}$. Then the following integration by parts formula holds: for every $\phi\in C_{b}^{\infty }({\mathbb{R}}^{d};\R)$, $Y\in \mathcal{S}$, $q\in\N$ and for every $\beta\in\{1,\ldots,d\}^q$  one has
\begin{equation*}
{\mathbb{E}}(\partial _{\beta}\phi(X)\,Y)={\mathbb{E}}(\phi(X)H_{\beta}^{q}(X,Y))
\end{equation*}%
where $\partial_\beta\phi(x)=\partial_{x^{\beta_q}}\cdots\partial_{x^{\beta_1}}\phi(x)$ and the weights $H^q_\beta(X,Y)$ are recursively given by:

\smallskip

$\bullet$
if $q=1$, then
$$
H^1_{\beta}(X,Y) \equiv H_{\beta}(X,Y)=\sum_{r=1}^{d}\Big(Y\gamma
_{X}^{r,\beta}LX^{r}-\sum_{k=1}^{n}\sum_{i=1}^{N}D_{(k,i)}(Y\gamma _{X}^{r,\beta})D_{(k,i)}X^{r}\Big),\quad\beta=1,\ldots,d;
$$
$\bullet$
if $q>1$, then
\begin{equation*}
H_{\beta}^{q}(X,Y)=H_{\beta_{q}}\big(%
X,H_{(\beta_{1},\ldots ,\beta_{q-1})}^{q-1}(X,Y)\big),\quad \beta\in\{1,\ldots,d\}^q.
\end{equation*}
Moreover, the following estimate holds: for every $\beta\in\{1,\ldots,d\}^q$ and $m\in\N$ then
\begin{equation}\label{gamma11}
\begin{array}{l}
\displaystyle
| H_{\beta}^{q}(X,Y)| _{m}\leq C\mathrm{\bf A}_{m+q}(X)^{q}|Y|_{m+q},
\quad\mbox{where}\smallskip\\
\mathrm{\bf A}_{l}(X)=\big(1\vee (\det \sigma_X)^{-1}\big)^{l+1}\big(1+|X|_{1,l+1}^{2d(l+2)}+|LX|_{l-1}^{2}\big),
\end{array}
\end{equation}%
$|\cdot|_m$ being defined in (\ref{abs11}).
\end{proposition}

We come now back to $S_{n}$, which we write as
$$
S_n=\frac 1{\sqrt n}\sum_{k=1}^n\big(\chi_kV_k+(1-\chi_k)W_k\big).
$$
For every $k=1,\ldots,n$ and $l,i=1,\ldots,N$, we have
\begin{equation*}
D_{(k,i)}S_{n}^{l}=\frac{1}{\sqrt{n}}\chi _{k}\I_{l=i}.
\end{equation*}%
As a consequence, we obtain
\begin{align}
&\sigma _{S_{n}}=\frac{1}{n}\sum_{k=1}^{n}\chi _{k}I,  \label{abs15}
\end{align}%
where $I$ denotes the identity matrix, and
\begin{align}
&\lambda _{S_{n}}= \frac{1}{n}%
\sum_{k=1}^{n}\chi _{k}  \label{abs16}
\end{align}%

The derivatives of order higher than two of $S_{n}$ are null, so we obtain
for every $q\in\N$

\begin{equation}
\left\vert S_{n}\right\vert _{1,q}^2\leq \frac{1}{%
n}\sum_{k=1}^{n}\chi _{k}\leq 1,\qquad \left\vert S_{n}\right\vert _{q}^2\leq
\left\vert S_{n}\right\vert^2 +\frac{1}{n}%
\sum_{k=1}^{n}\chi _{k}\leq \left\vert S_{n}\right\vert^2+ 1,  \label{abs17}
\end{equation}%
and consequently%
\begin{equation}
\left\Vert S_{n}\right\Vert _{1,q,p}\leq 1,\qquad
\left\Vert S_{n}\right\Vert _{q,p}\leq \left\Vert S_{n}\right\Vert _{p}+1.  \label{abs18}
\end{equation}
In particular, $\left\Vert S_{n}\right\Vert _{1,q,p}$ is finite for every $q,p$ whereas $\left\Vert S_{n}\right\Vert _{q,p}$ is finite according to $F\in L^p(\Omega)$.

Let us now compute $LS_{n}.$ We have%
\begin{eqnarray*}
-LS_{n}^{l}
&=&\sum_{k=1}^{n}\sum_{i=1}^{N}D_{(k,i)}D_{(k,i)}S_{n}^{l}+\sum_{k=1}^{n}%
\sum_{i=1}^{N}D_{(k,i)}S_{n}^{l}\partial _{i}\ln \psi _{r_0/2}(\left\vert
V_{k}-v_{0}\right\vert ) \\
&=&\frac{1}{\sqrt{n}}\sum_{k=1}^{n}\chi _{k}\partial _{l}\ln \psi _{r_0/2}(\left\vert V_{k}-v_{0}\right\vert ).
\end{eqnarray*}%
We now estimate $\left\Vert LS_{n}\right\Vert _{q,p}.$

\begin{lemma}
For every $q\in \N$, there exists a universal constant $C_{q}$ such that%
\begin{equation}
\left\Vert LS_{n}\right\Vert _{q,p}\leq \frac{C_{q}}{r_0^{q+1}}.  \label{abs19}
\end{equation}
\end{lemma}

\textbf{Proof}. The basic fact in our calculus is that
\begin{align*}
\E(\partial _{i}\ln \psi _{r_0/2}(V_{k}-v_{0}))
=&\frac{\varepsilon_0 }{m_0}%
\int_{\R^{N}}\partial _{i}\ln \psi _{r_0/2}(\left\vert v-v_{0}\right\vert
)\times \psi _{r_0/2}(\left\vert v-v_{0}\right\vert )dv  \\ 
=&\frac{\varepsilon_0 }{m_0}\int_{\R^{N}}\partial _{i}\psi _{r_0/2}(\left\vert
v-v_{0}\right\vert )dv=0.
\end{align*}%
We denote%
\begin{equation*}
Q_{k}=\nabla \ln \psi _{r_0/2}(V_{k}-v_{0})
\end{equation*}%
and we have
\begin{equation*}
\E(Q_{k}^{l})=\E(\partial _{l}\ln \psi _{r_0/2}(\left\vert V_{k}-v_{0}\right\vert ))=0.
\end{equation*}%
So $\sum_{k=1}^n\chi_kQ_k^l$, $n\in\N$, is a martingale and the Burkholder's inequality gives
\begin{equation*}
\E(| LS_{n}^{l}|^{p})
=\E\Big(\Big\vert \frac{1}{\sqrt{n}}%
\sum_{k=1}^{n}\chi _{k}Q_{k}^{l}\Big\vert ^{p}\Big)
\leq C\E\Big(\Big(\frac{1}{n}%
\sum_{k=1}^{n}\chi _{k}\left\vert Q_{k}^{l}\right\vert ^{2}\Big)^{p/2}\Big)
\leq
\frac{C}{n}\sum_{k=1}^{n}\E(|Q_{k}^{l}|^{p}).
\end{equation*}%
By (\ref{abs3})%
\begin{equation*}
\E(\vert Q_{k}^{l}\vert ^{p})\leq C\,\frac{%
1}{r_0^{p}}
\end{equation*}%
so that
\begin{equation*}
\left\Vert LS_{n}\right\Vert _{p}\leq \frac{C}{r_0}.
\end{equation*}%
We go further and we compute $D_{(k,i)}LS_{n}.$ We have
$$
-D_{(k,i)}LS_{n}^l =\frac{1}{\sqrt{n}}\sum_{k^{\prime }=1}^{n}\chi
_{k^{\prime }}D_{(k,i)}\partial_l\ln \psi _{r_0/2}(\left\vert
V_{k^{\prime }}-v_{0}\right\vert )
=\frac{1}{\sqrt{n}}\chi _{k}D_{(k,i)}\partial_l \ln \psi _{r_0/2}(\left\vert V_{k}-v_{0}\right\vert )
$$
so that
\begin{eqnarray*}
\left\vert DLS_{n}\right\vert _{1}^{2} &\leq &\left\vert LS_{n}\right\vert
^{2}+\frac{1}{n}\sum_{k=1}^{n}\sum_{i=1}^{N}\left\vert D_{(k,i)}\nabla \ln \psi _{r_0/2}(\left\vert V_{k}-v_{0}\right\vert )\right\vert^2  \\
&\leq &\left\vert LS_{n}\right\vert ^{2}+\frac{C%
}{n}\sum_{k=1}^{n}\sum_{i,j=1}^{N}\left\vert \partial _{i}\partial _{j}\ln
\psi _{r_0/2}(\left\vert V_{k}-v_{0}\right\vert )\right\vert^2 .
\end{eqnarray*}%
Once again using (\ref{abs3}) we obtain
\begin{equation*}
\left\Vert \partial _{i}\partial _{j}\ln \psi _{r_0/2}(\left\vert
V_{k}-v_{0}\right\vert )\right\Vert _{p}\leq \frac{C}{r_0^{2}}
\end{equation*}%
and consequently%
\begin{equation*}
\left\Vert LS_{n}\right\Vert _{1,p}\leq \frac{C}{r_0^{2}}.
\end{equation*}%
For higher order norms, the estimates are similar. $\square $

\medskip

We add a final property on the behavior of the Malliavin covariance matrix that will be used in next Section \ref{sect-regularization}.

\begin{lemma}\label{est-sigma}
Suppose that $\mu _{F}\succeq \leb_{N}.$
There exists a
universal constant $C$ such that for every $n\in \N$ and every
\begin{equation}
\varepsilon \leq \varepsilon _{\ast }=2^{-N}m_0^{N}
\label{abs28}
\end{equation}%
then
\begin{equation}
\P(\det \sigma _{S_{n}}\leq \varepsilon )\leq C\exp (-\frac{n}{4(\frac{1}{m_0}
-1)}),  \label{abs29}
\end{equation}%
$m_0$ being defined in (\ref{abs5}).
\end{lemma}

\textbf{Proof.}
Using (\ref{abs16})%
$$\P(\det \sigma _{S_{n}} \leq \varepsilon )
\leq
\P(\lambda _{S_{n}}\leq \varepsilon ^{1/N})
= \P\Big(\frac{1}{n}\sum_{k=1}^{n}\chi _{k}\leq \varepsilon ^{1/N}\Big)
=\P\Big(\frac{1}{n}\sum_{k=1}^{n}(\chi _{k}-m_0)\leq \varepsilon ^{1/N}-m_0\Big).
$$
Since $\varepsilon ^{1/N}\leq \frac{1}{2}m_0$, the
above term is upper bounded by%
\begin{equation*}
\P\Big(\frac{1}{n}\sum_{k=1}^{n}(\chi _{k}-m_0)\leq -\frac{1}{2}m_0\Big)
=\P\Big(\frac{1}{n^{1/2}}\sum_{k=1}^{n}\frac{\chi _{k}-m_0}{v_{m_0}}\leq -n^{1/2}\frac{m}{2v_{m_0}}\Big)
\end{equation*}%
with $v_{m_0}=(m_0(1-m_0))^{1/2}=\Var(\chi _{k}).$ We denote by $a=n^{1/2}\frac{m_0}{%
2v_{m_0}}$ and we use the Berry--Esseen theorem in order to upper bound this
quantity by%
\begin{equation*}
C\int_{-\infty }^{a}\exp \big(-x^{2}/{2}\big)dx
\leq C^{\prime }\exp \Big(-\frac{a^{2}}{4}\Big)
=C^{\prime }\exp \Big(-\frac{n}{4(\frac{1}{m_0}-1)}\Big).
\end{equation*}
$\square$

\subsection{Proof of Theorem \ref{TV}}

We need now a localized variant of Lemma 2.5 and Theorem 2.7 in \cite{[BC]}. So, we start with the basic definitions.

We consider a localizing r.v. $\Theta$  taking values in $[0,1]$ of the form
\begin{equation}\label{Theta}
\Theta=\psi_a(Z),\quad a>0,\quad Z\in\SS,
\end{equation}
$\psi_a$ being defined in (\ref{abs2}).  We set $\P_\Theta$ and $\E_\Theta$ through
$$
d\P_\Theta=\Theta d\P\quad\mbox{and}\quad \E_\Theta=\mbox{expectation w.r.t. }\P_\Theta.
$$
For $X\in\SS^d$, we define the localized Sobolev norms
$$
\|X\|_{p,\Theta}=\E_\Theta(|X|^p)^{1/p},\quad
\|X\|_{1,m,p,\Theta}=\E_\Theta(|X|_{1,m}^p)^{1/p}
\quad\mbox{and}\quad
\|X\|_{m,p,\Theta}=\E_\Theta(|X|_m^p)^{1/p},
$$
$|X|_{1,m}$ and $|X|_m$ being given in (\ref{abs11}), and we set
\begin{equation}\label{A}
A_{p,\Theta}(X)=\|X\|_{3,p,\Theta}+\|LX\|_{1,p,\Theta}.
\end{equation}
We also consider the law of a $d$-dimensional r.v. $X$ under $\P_\Theta$: it is the measure in $\R^d$ defined as
$$
\mu_{X,\Theta}(dx)=\P_{\Theta}(X\in dx).
$$
We allow the case $a=+\infty$ in (\ref{Theta}): this gives $\Theta\equiv 1$, so $\P_\Theta\equiv \P$ and no localization is taken into account.

Finally, for $k\in\N$, we define the distance $d_k$ between two measures $\mu,\nu$ in $\R^d$  as
\begin{equation}\label{dk}
d_{k}(\mu,\nu)=\sup \Big\{\Big|\int fd\mu-\int fd\nu\Big|\,:\,\|f\|_{k,\infty}\leq 1\Big\}
\end{equation}
where $\|f\|_{k,\infty}=\sum_{0\leq |\alpha|\leq k}\|\partial_\alpha f\|_{\infty}$. Then we have $d_0=d_{TV}$ and $d_1=d_{FM}$ (Fortet--Mourier distance).

The following result is a localized version of Lemma 2.5 in \cite{[BC]}.
Here, $\gamma _{\delta }$ denotes the density of the centred
normal law of covariance $\delta \times I$ on ${\mathbb{R}}^{d}$, $\delta >0$ ($I$ denoting the identity matrix) and $f\ast\gamma_\delta$ denotes the convolution between $f$ and $\gamma_\delta$.

\begin{lemma}
\label{Lemma2.5}
Let $\Theta$ be a localizing r.v. as in (\ref{Theta}). Then,
for every $\varepsilon >0,\delta >0$,  $X\in \SS^d$ and for every bounded and measurable function $f:{\mathbb{R}}^{d}\rightarrow {\mathbb{R}}$ one has
\begin{equation}
\left\vert {\mathbb{E}_\Theta}(f(X))-{\mathbb{E}_\Theta}(f\ast \gamma _{\delta
}(X))\right\vert \leq C\left\Vert f\right\Vert _{\infty }\Big({\mathbb{P}}_\Theta%
(\sigma _{X}<\varepsilon )+\frac{\sqrt{\delta }}{\varepsilon ^{p}}%
(1+A_{p,\Theta}(X))^{a}\Big)
\label{Main5}
\end{equation}%
where $A_{p,\Theta}(X)$ is defined in (\ref{A}) and $C,p,a>0$ are suitable universal constants depending on the dimension $d$ only.
\end{lemma}

\textbf{Proof.} The proof is identical to the one of Lemma 2.5 in \cite{[BC]} (the case $\Theta\equiv 1$ being the same result): just consider the localized measure $\P_\Theta$ instead of $\P$ in the proof of Lemma 2.5 in \cite{[BC]} (namely, replace the localizing variable $\Theta_\varepsilon$ therein with $\Theta_\varepsilon\Theta$).
$\square$

\medskip

We state now a variant of Theorem 2.7 in \cite{[BC]} that takes into account  localizations.

\begin{theorem}\label{Th2.7}
Let $\Theta,U$ be localizing r.v.'s as in (\ref{Theta}) and let $X, Y\in \SS^d$
be such that $A_{l,\Theta}(X)$, $A_{l,U}(Y)<\infty$, such quantities being defined in (\ref{A}). Let $\mu_{X,\Theta}$ denote the law of $X$ under $\P_\Theta$ and let $\mu_{Y,U}$ denote the law of $Y$ under $\P_U$. Let $k\in\N$. Then there exist some universal constants $C,p,a,b>0$ (independent of $\Theta$, $U$, $X$, $Y$, $k$) such that
\begin{equation}  \label{New1}
\begin{array}{rcl}
d_{0}(\mu_{X,\Theta},\mu_{Y,U}) & \leq & \displaystyle\frac{C}{\varepsilon ^{b}}%
\big(1+A_{l,\Theta}(X)+A_{l,U}(Y)\big)^{a}\big(d_{k}(\mu_{X,\Theta},\mu_{Y,U})%
\big)^{\frac{1}{k+1}}+ \smallskip \\
&  & +C{\mathbb{P}_\Theta}(\det \sigma _{X} <\varepsilon )+C{\mathbb{P}_U}(\det
\sigma _{Y}<\varepsilon ),
\end{array}%
\end{equation}
\end{theorem}

\textbf{Proof.}
We take a bounded and measurable function $f$ and we write
\begin{align*}
|\E_\Theta(f(X))-\E_U(f(Y))|
\leq&
|\E_\Theta(f(X))-\E_\Theta(f\ast\gamma_\delta(X))|
+|\E_U(f(X))-\E_U(f\ast\gamma_\delta(Y))|\\
&+|\E_\Theta(f\ast\gamma_\delta(X))-\E_U(f\ast\gamma_\delta(Y))|\\
=:
&I_\Theta(X)+I_U(Y)+I_{\Theta,U}(X,Y).
\end{align*}
By using (\ref{Main5}) we get
\begin{align*}
I_\Theta(X)
+I_U(Y)
\leq &C\left\Vert f\right\Vert _{\infty }\Big({\mathbb{P}}_\Theta%
(\sigma _{X}<\varepsilon )+{\mathbb{P}}_U
(\sigma _{Y}<\varepsilon )+\frac{\sqrt{\delta }}{\varepsilon ^{p}}%
(1+A_{p,\Theta}(X)+A_{p,U}(Y))^{a}\Big)
\end{align*}
Moreover, by recalling that $\|f\ast \gamma
_{\delta }\|_{k,\infty }\leq C\delta ^{-k/2}\|f\|_\infty$, we have
\begin{align*}
I_{\Theta,U}(X,Y)
\leq C\delta ^{-k/2}\|f\|_\infty d_k(\mu_{X,\Theta}, \mu_{Y,U}).
\end{align*}
Following the proof of Theorem 2.7 in \cite{[BC]}, we now insert everything, optimize w.r.t. $\delta$  and we get the result.
$\square$

\begin{remark}\label{rmGaussian}
Lemma \ref{Lemma2.5} and Theorem \ref{Th2.7} are valid not only with the basic noise $V_1,\ldots,V_n$ introduced in Section \ref{AMC}. Actually, both results remains true whenever the basic noise fulfils the abstract integration by parts framework developed in Section 2.1 of \cite{[BC]}, the one considered in this paper being a particular case.
\end{remark}

We are finally ready for the following.

\medskip

\textbf{Proof of Theorem \ref{TV}}.
Let $G$ denote a standard normal r.v. in $\R^N$. For each $K\geq 1$ set
$$
\Theta _{n,K}=\psi _{K}(S_{n}),\quad d\P_{\Theta _{n,K}}=\Theta
_{n,K}d\P\quad\mbox{and}\quad
\Theta _{K}=\psi _{K}(G),\quad d\P_{\Theta _{K}}=\Theta
_{K}d\P,
$$
$\psi_K$ being defined in (\ref{abs2}). Let $\mu _{n,K}$ be the law of $S_n$ under $\P_{\Theta_{n,K}}$ and $\mu _{K}$ be the law of $G$ under $\P_{\Theta_{K}}$, that is,
$$
\mu _{n,K}(dx) =\P_{\Theta _{n,K}}(S_{n}\in dx)\quad\mbox{and}\quad
\mu _{K}(dx) = \P_{\Theta _{K}}(G\in dx)
$$
Consider a measurable function $f:\R^{N}\rightarrow \R$
such that $\left\Vert f\right\Vert _{\infty }\leq 1$. We write
\begin{eqnarray*}
\left\vert \E(f(S_{n}))-\E(f(G))\right\vert  &\leq &\left\vert
\E(f(S_{n})(1-\Theta _{n,K}))\right\vert +\left\vert
\E(f(G)(1-\Theta _{K}))\right\vert  \\
&&+\left\vert \E(f(S_{n})\Theta _{n,K})-\E(f(G)\Theta _{K})\right\vert .
\end{eqnarray*}%
Using the Chebyshev's inequality,
\begin{equation*}
\left\vert \E(f(S_{n})(1-\Theta _{n,K}))\right\vert \leq
\left\Vert f\right\Vert _{\infty }\P(\left\vert S_{n}\right\vert \geq 2K)\leq
\frac{C}{K^{2}}\left\Vert f\right\Vert _{\infty }
\end{equation*}%
and a similar estimates holds for $|\E(f(G)(1-\Theta _{K}))|$. We conclude that%
\begin{equation*}
\sup_{\left\Vert f\right\Vert _{\infty }\leq 1}\left\vert
\E(f(S_{n}))-\E(f(G))\right\vert \leq \frac{C}{K^{2}}+\sup_{\left\Vert
f\right\Vert _{\infty }\leq 1}\left\vert \E(f(S_{n})\Theta
_{n,K})-\E(f(G)\Theta _{K})\right\vert .
\end{equation*}%
We obtain%
\begin{equation*}
\limsup_{n\rightarrow \infty }\sup_{\left\Vert f\right\Vert _{\infty
}\leq 1}\left\vert \E(f(S_{n}))-\E(f(G))\right\vert \leq \frac{C}{K^{2}}
+\limsup_{n\rightarrow \infty }d_{TV}(\mu_{n,K},\mu_K)
\end{equation*}%
for every $K\geq 1$. If we show that, for each fixed $K$, $d_{TV}(\mu_{n,K},\mu_K)\to 0$ as $n\to\infty$, the statement will follow by letting $K$ go to $+\infty$. So, we study $d_{TV}(\mu_{n,K},\mu_K)$, for a fixed $K>1$.

We use Theorem \ref{Th2.7} with $\Theta=\Theta_{n,K}$, $X=S_n$, $U=\Theta_K$ and $Y=G$. Here, the noise includes the Gaussian r.v. $G$, so we add it to the underlying noise (recall Remark \ref{rmGaussian}) in a standard way -- we stress this trick because it will be used also in the sequel, for example, in Lemma \ref{4.4-new}.

Without loss of generality, we assume that $G$ is defined on the same probability space and is independent of $V_1,\ldots,V_n$. We consider as basic noise the one coming from $(G,V_1,\ldots,V_n)$. For $X=\phi(G,V_1,\ldots,V_n)$ with $\phi\in C_b^{\infty}(\R^{N(1+n)};\R)$, we set
$$
D_{(0,i)}X=\frac{\partial}{\partial G^i}\phi(G,V_1,\ldots,V_n)
$$
and $D_{(k,i)}$ for $k=1,\ldots,n$ as in (\ref{abs9}). The Ornstein--Uhlenbeck generator takes into account the contribution from the standard Gaussian $G$, so it becomes
\begin{align*}
-LX
=&\sum_{i=1}^{N}D_{(0,i)}D_{(0,i)}X-\sum_{i=1}^{N}D_{(0,i)}X \, G^i\\
&+\sum_{k=1}^{n}\sum_{i=1}^{N}D_{(k,i)}D_{(k,i)}X+\sum_{k=1}^{n}%
\sum_{i=1}^{N}D_{(k,i)}X\partial _{i}\ln \psi _{r_0/2}(\left\vert
V_{k}-v_{0}\right\vert ).
\end{align*}
And if $X$ is a random vector in $\R^d$, the associated Malliavin covariance matrix is
$$
\sigma _{X}^{i,j}
=\sum_{k=0}^{n}\sum_{r=1}^{N}D_{(k,r)}X^{i}\times D_{(k,r)}X^{j},\quad i,j=1,\ldots,d.
$$
It is standard to see that the above quantities bring to an abstract Malliavin calculus as developed in \cite{[BC]}. Of course, when the randomness does not depend on $G$ then everything agrees with what developed in Section \ref{AMC} and when the randomness does not depend on $V$ then we get the standard Gaussian--Malliavin calculus. So, we use Remark \ref{rmGaussian} and we apply Theorem \ref{Th2.7}. In order to use (\ref{New1}), we need to study $A_{\Theta_{n,K}}(S_n)$ and $A_{\Theta_K}(G)$. By  (\ref{abs17}) and by recalling that $1_{\{\Theta _{n,K}\neq 0\}}|S_{n}|\leq 2K$, we obtain
$$
\left\Vert S_{n}\right\Vert _{q,p,\Theta _{n,K}}+\left\Vert
LS_{n}\right\Vert _{q-2,p,\Theta _{n,K}}\leq CK.
$$
Standard computations give $\|G\|_{q,p, \Theta_K}+\|LG\|_{q-2,p, \Theta_K}
\leq \|G\|_{q,p}+\|LG\|_{q-2,p}\leq C$, so we can write
$$
A_{\Theta_{n,K}}(S_n) + A_{\Theta_K}(G)\leq CK,
$$
$C>0$ being independent of $K$ and $n$. Moreover, $\sigma_G$ is the identity matrix. And
since $\left\vert \Theta
_{n,K}\right\vert \leq 1$ and $\chi _{k},k\in \N$ are i.i.d, the law of
large numbers says that for $\varepsilon^{1/N} <\E(\chi _{k})=m_0$ one has%
\begin{equation*}
\limsup_{n\to\infty}\P_{\Theta _{n,K}}(\det \sigma _{S_{n}}\leq \varepsilon )
\leq \limsup_{n\to\infty}\P_{\Theta _{n,K}}(\lambda _{S_{n}}\leq \varepsilon^{1/N} )\leq
\limsup_{n\to\infty}\P\Big(\frac{1}{n}\sum_{k=1}^{n}\chi _{k}\leq \varepsilon^{1/N}\Big)=0
\end{equation*}%
in which we have used (\ref{abs16}). We apply now Theorem \ref{Th2.7} with $k=1$ and $\varepsilon <1\wedge m_0^N$: by passing to the limit in (\ref{New1}) we obtain
$$
\limsup_{n\to\infty}d_{0}(\mu_{n,K},\mu_{K})
\leq \displaystyle\frac{C}{\varepsilon ^{a}}%
\big(1+CK\big)^{b}\limsup_{n\to\infty}d_{FM}(\mu_{n,K},\mu_{K})^{1/2}.
$$
So, it remains to show that $d_{FM}(\mu_{n,K},\mu_{K})\to 0$ as $n\to\infty$.
Since $\psi _{K}\in C_{c}(\R^{N})$, the CLT gives
\begin{equation*}
\lim_{n}\E_{\Theta_{n,K}}(f(S_{n}))
=\lim_{n}\E(\psi _{K}(S_{n})f(S_{n}))=\E(\psi _{K}(G)f(G))
=\E_{\Theta_K}(f(G))
\end{equation*}%
for every $f\in C(\R^{d})$. So, if we define the probability laws
$$
\hat\mu_{n,K}(dx)=\frac 1{\E(\Theta_{n,K})}\mu_{n,K}(dx)\quad\mbox{and}\quad
\hat\mu_{K}(dx)=\frac 1{\E(\Theta_{K})}\mu_{K}(dx),
$$
we get $\hat\mu_{n,K}\to \hat\mu_{K}$ weakly as $n\to\infty$. Since weak convergence of probability laws is equivalent to convergence in $d_{FM}$ (see e.g. Theorem 11.3.3 in \cite{dudley}), we have
$d_{FM}(\hat\mu_{n,K},\hat\mu_{K})\to 0$ as $n\to\infty$. Finally, straightforward computations give
$$
d_{FM}(\mu_{n,K},\mu_{K})\leq \big|\E(\Theta_{n,K})-\E(\Theta_K)\big|+d_{FM}(\hat\mu_{n,K},\hat\mu_{K})\to 0
$$
as $n\to\infty$, and the statement follows.
$\square $

\begin{remark}\label{noId1}
We note that if $C(F)$ was not the identity matrix then (\ref{abs15}) and (\ref{abs16}) would become
$$
\sigma_{S_n}=\frac 1n\sum_{k=1}^n\chi_k \widehat{C}(F)\quad\mbox{and}\quad \lambda_{S_n}= \underline{\lambda}(F)\frac 1n\sum_{k=1}^n\chi_k
$$
respectively, where $\widehat{C}(F)=C(F)^{-1}$ and $\underline{\lambda}(F)$ is the smallest eigenvalue of $\widehat{C}(F)$. This means that the estimates in (\ref{abs18}) and (\ref{abs19}) continue to hold up to a multiplying constant that now depends on $\underline{\lambda}(F)$ and $\overline{\lambda}(F)$ as well, the latter denoting the largest eigenvalue of $\widehat{C}(F)$.
\end{remark}

\section{Asymptotic expansion}\label{sect-asympt}

The aim of this section is to prove Theorem \ref{Speed-r} and  Theorem  \ref{main-th-1}. We first study the case of smooth functions and then, using a regularizing argument, we will be able to deal with general functions.

\subsection{The development for smooth test functions}\label{sect-dev-smooth}

We recall that we are assuming that the r.v. $F$ has null mean and non-degenerate covariance matrix, that we have set equal to the identity matrix. And we have set
\begin{equation*}
F_{i} =\chi_iV_i+(1-\chi_i)W_i
\end{equation*}%
so that $S_{n}=\frac{1}{\sqrt{n}}\sum_{i=1}^{n}F_i=\frac 1{\sqrt n}\sum_{i=1}^n(\chi_iV_i+(1-\chi_i)W_i)$.
Moreover we consider $%
G_{i}=(G_{i}^{1},...,G_{i}^{N}),i\in \N$, some independent standard normal random variables in $\R^N$. For $k\in \{0,1,...,n\}$, we define
\begin{equation}
S_{n}^{k}=\frac{1}{\sqrt{n}}\Big(\sum_{i=1}^{k}F_i+\sum_{i=k+1}^{n}G_{i}\Big),%
\qquad \widehat{S}_{n}^{k}=\frac{1}{\sqrt{n}}\Big(\sum_{i=1}^{k-1}F_i+%
\sum_{i=k+1}^{n}G_{i}\Big)  \label{abs24}
\end{equation}%
in which we use the convention that the sums are null when done on the indexes $i\in\{i_0,\ldots,i_1\}$ with $i_0>i_1$. Therefore, one has
\begin{equation*}
S_{n}^{n}=S_{n}\quad\mbox{and}\quad S_{n}^{0}=\frac{1}{\sqrt{n}}\sum_{i=1}^{n}G_{i}
\end{equation*}%
and $S_{n}^{0}$ is a standard normal random variable in $\R^N$. Moreover,
\begin{equation}
S_{n}^{k}=\widehat{S}_{n}^{k}+\frac{F_k}{\sqrt n}\quad\mbox{and}\quad S_{n}^{k-1}=\widehat{S}%
_{n}^{k}+\frac{G_{k}}{\sqrt n}.  \label{abs24'}
\end{equation}%

In the sequel, we will use the following notation. For a multiindex $\alpha
=(\alpha _{1},...,\alpha _{r})\in \{1,...,N\}^{r}$ and $x=(x^{1},...,x^{N})$
we denote $x^{\alpha }=\prod_{i=1}^{r}x^{\alpha _{i}}$. We also denote by $\partial _{\alpha }=\partial _{x^{\alpha _{1}}}...\partial _{x^{\alpha
_{r}}} $ the derivative corresponding to $\alpha $ and by $\left\vert
\alpha \right\vert =r$ the length of $\alpha .$ We allow $\alpha$ to be the null multiindex: in this case, we set $|\alpha|=0$, $\partial_\alpha f=f$ and $x^\alpha=1$.

Moreover, we will use the
following form of the Taylor formula of order $r\in \N:$ for $f\in C^{r+1}(\R^N)$,
\begin{equation}
f(x+y)=f(x)+\sum_{p=1}^{r}\frac{1}{p!}\sum_{\left\vert \alpha \right\vert
=p}\partial _{\alpha }f(x)y^{\alpha }+U_{r}f(x,y)  \label{Taylor1}
\end{equation}%
with%

\begin{equation}
U_{r}f(x,y)
=\frac 1{r!}\sum_{\left\vert \alpha \right\vert =r+1}y^{\alpha}\int_0^1(1-\lambda)^r\partial _{\alpha }f(x+\lambda y)d\lambda
\label{Taylor2}
\end{equation}

We notice that for some $c_r>0$ it holds
\begin{equation}
|U_{r}f(x,y)|\leq c_r |y|^{r+1}\|f\|_{r+1,\infty},
\label{Taylor2'}
\end{equation}
where $\|\cdot \|_{r+1,\infty}$ is the usual norm on $C^{r+1}_b(\R^N)$: $\|f\|_{r+1,\infty}=\sum_{|\alpha|\leq r+1}\|\partial_\alpha f\|_\infty$.

For a multiindex $\alpha =(\alpha _{1},...,\alpha _{r})\in \{1,...,N\}^{r}$, that is, $|\alpha|=r$, we now set
\begin{equation}\label{DO1}
\begin{array}{l}
\displaystyle
\Delta _{\alpha }=\E(F^{\alpha })-\E(G^{\alpha })=\E\big(\prod_{i=1}^{r}F^{\alpha
_{i}}\big)-\E\big(\prod_{i=1}^{r}G^{\alpha _{i}}\big),\smallskip\\
\displaystyle
\mbox{$\theta _{\alpha}=1$ if $r$ is even and $\alpha_{2j-1}=\alpha_{2j}$ for every $j=1,\ldots r/2$, otherwise $\theta _{\alpha}=0$.}
\end{array}
\end{equation}%
For $\alpha=\emptyset$, that is, $r=0$, we set $\Delta_\emptyset=0$ and $\theta_\emptyset=1$. It is clear that $\Delta_\alpha=0$ for $|\alpha|\leq 2$ and, for $r\geq 3$, the assumption $\sup_{|\alpha|\leq r}|\Delta_\alpha|=0$ means that all moments of $F$  up to order $r$ (and not only up to order 2) agree with the moments of a standard Gaussian random variable.

We now introduce the basic differential operators which appear in the
asymptotic expansion: we set
\begin{equation}\label{Psi}
\Psi_t=\sum_{p=0}^ t
%
%
\frac  {(-1)^{\frac{t-p}2}}{2^{\frac{t-p}2}p!(\frac{t-p}2)!}\sum_{|\alpha|=p}
\sum_{|\beta|={t-p}} \theta_\beta\Delta_\alpha\partial_\beta\partial_\alpha,\quad t=0, 1,2,\ldots.
\end{equation}
%
%
Recall that $\theta_\beta$ is null when $t-p$ is odd, so the sum actually runs on the indexes $p$ such that $\frac{t-p}2\in\N$.
The property $\Delta_\alpha=0$ if $|\alpha|\leq 2$ gives that the sum in (\ref{Psi}) actually starts from $p=3$, so we have
$$
\Psi_t=0\quad\mbox{if}\quad t=0,1,2
\quad\mbox{and}\quad
\Psi_t=\sum_{p=3}^t\frac  {(-1)^{\frac{t-p}2}}{2^{\frac{t-p}2}p!(\frac{t-p}2)!}\sum_{|\alpha|=p}
\sum_{|\beta|={t-p}} \theta_\beta\Delta_\alpha\partial_\beta\partial_\alpha\quad\mbox{if}\quad t\geq 3.
$$
From now on, we use the convention $\sum_{p=3}^t(\cdot)=0$ if $t<3$. So, for example we can write
$$
\Psi_t=\sum_{p=3}^t\frac  {(-1)^{\frac{t-p}2}}{2^{\frac{t-p}2}p!(\frac{t-p}2)!}\sum_{|\alpha|=p}
\sum_{|\beta|={t-p}} \theta_\beta\Delta_\alpha\partial_\beta\partial_\alpha,\quad t=0, 1,2,\ldots
$$
We note that $\Psi_t=0$ for all $t$ when $\Delta_\alpha=0$ for all $\alpha$, that is when all the moments of $F$ agree with the moments of the standard Gaussian law. And moreover, for every $t\geq 3$ and $q\geq 0$ there exists $C_{t,q}>0$ such that if $f\in C^{t+q}_b$ then
\begin{equation}\label{est-Psit}
\|\Psi_t f\|_{q,\infty}
\leq
C_{t,q}\sup_{|\alpha|\leq t}|\Delta_\alpha|\times\|f\|_{t+q,\infty}.
\end{equation}
We also define the following objects (``remainders''): for $r\in\N$ and $f\in C_b^{r+2}(\R^N)$,
\begin{equation}\label{Rkf}
\begin{array}{l}
\displaystyle
\mathcal{R}_{r,n}^kf
=n^{\frac{r+1}2}\Big[\E\Big(U_{r}f\Big(\widehat S^k_n,\frac{F_k}{\sqrt n}\Big)\Big)
-\E\Big(U_{r}f\Big(\widehat S^k_n,\frac{G_k}{\sqrt n}\Big)\Big)\Big]\smallskip\\
\displaystyle
+\sum_{p=0}^r
n^{-([\frac{r-p}2]+\frac 12-\frac{r-p}2)}
\times
\frac {(-1)^{[\frac{r-p}2]+1}}{  p![\frac{r-p}2]! 2^{[\frac{r-p}2]+1}}
\!\!\!\!\sum_{\mbox{\scriptsize{$\begin{array}{c}|\alpha|=p\\|\beta|={2[\frac{r-p}2]+2}\end{array}$}}}\!\!\!\!\!\!\!\!
\Delta_\alpha
\theta_\beta\int_0^1 s^{[\frac{r-p}2]}
\E\Big(\partial_\beta\partial_\alpha f\Big(\widehat S^{k}_n+\sqrt s\,\frac {G_k}{\sqrt n}\Big)\Big)ds,
\end{array}
\end{equation}
$U_{r}f$ being defined in (\ref{Taylor2}). Note that the second term of the above right hand side is equal to zero if $r<3$. Moreover, $[\frac{r-p}2]+\frac 12-\frac{r-p}2\in\{0,\frac 12\}$, hence $n^{-([\frac{r-p}2]+\frac 12-\frac{r-p}2)}\leq 1$.

\begin{remark}\label{rem-R01}
We note here if $F\in L^2$ then for every $f\in C^2_b$ one has
$$
\mathcal{R}_{0,n}^kf
= \frac 1{\sqrt n} \mathcal{R}_{1,n}^kf.
$$
And if $F\in L^3(\Omega)$ then for every $f\in C^3_b$ one has
\begin{equation}\label{R01}
\mathcal{R}_{0,n}^kf
= \frac 1{\sqrt n} \mathcal{R}_{1,n}^kf
= \frac 1n \mathcal{R}_{2,n}^kf.
\end{equation}
In fact, for every $r\geq 0$, if $f\in C^{r+2}_b$ then
$$
U_{r}f(x,y)=U_{r+1}f(x,y)-\frac 1{(r+1)!}\sum_{|\alpha|=r+1}y^\alpha\partial_\alpha f(x).
$$
Therefore, for $r=0$, $F\in L^2$ and $f\in C^2_b$ we obtain
$$
\mathcal{R}_{0,n}^kf
=\sqrt n\Big[\E\Big(U_{1}f\Big(\widehat S^k_n,\frac{F_k}{\sqrt n}\Big)\Big)
-\E\Big(U_{1}f\Big(\widehat S^k_n,\frac{G_k}{\sqrt n}\Big)\Big)\Big]-\sqrt n\sum_{|\alpha|=1}\E\Big(\Big[\Big(\frac{F_k}{\sqrt n}\Big)^\alpha
-\Big(\frac{G_k}{\sqrt n}\Big)^\alpha\Big]f(\widehat S^k_n)\Big).
$$
Since $\widehat S^k_n$ is independent of $F_k$ and $G_k$ and since $\Delta_\alpha=0$ for $|\alpha|=1$ we get
$\E([(F_k\big)^\alpha
-(G_k)^\alpha]f(\widehat S^k_n))
=\Delta_\alpha\E(f(\widehat S^k_n))
=0$,
 so that
$$
\mathcal{R}_{0,n}^kf
= \frac 1{\sqrt n} \mathcal{R}_{1,n}^kf.
$$
As for (\ref{R01}), one uses $\Delta_\alpha=0$ for $|\alpha|=2$ and the statement is proved similarly.
\end{remark}

Since $\E(f(S_n))-\E(f(G))=\E(f(S^n_n))-\E(f(S^0_n))$, we study $\E(f(S^{k}_n))-\E(f(S^{k-1}_n))$ for $k=1,\ldots,n$ and then apply a recurrence argument.

\begin{lemma}\label{development}
Let $n\in \N,1\leq k\leq n$ and $r\in\N$. If $F\in L^{r+1}(\Omega)$ then for every  $f\in C^{r+1}_b(\R^{N})$ one has
\begin{equation}
\E(f(S^{k}_n))-\E(f(S^{k-1}_n))
=\sum_{p=3}^r\frac 1{p! n^{p/2}}\sum_{|\alpha|=p}
\E\big(\partial_\alpha f(\widehat S^{k}_n)\big)\Delta_\alpha+\frac 1{n^{(r+1)/2}}
\widetilde{\mathcal{R}}^k_{r,n}f
 \label{taylor}
\end{equation}%
where
$$
\widetilde{\mathcal{R}}^k_{r,n}f
=n^{\frac{r+1}2}\Big[\E\Big(U_{r}f\Big(\widehat S^k_n,\frac{F_k}{\sqrt n}\Big)\Big)
-\E\Big(U_{r}f\Big(\widehat S^k_n,\frac{G_k}{\sqrt n}\Big)\Big)\Big].
$$
\end{lemma}

\textbf{Proof}. We will use the Taylor formula (\ref{Taylor1}). Since $S_n^k=\widehat S_n^k+\frac{F_k}{n^{1/2}}$ and $F_k$ is independent of $\widehat{S}_{n}^{k}$, we obtain%
\begin{equation*}
\E(f(S_{n}^{k}))=\E(f(\widehat{S}_{n}^{k}))
+\sum_{p=1}^{r}\frac{1}{p!n^{p/2}}\sum_{\left\vert \alpha \right\vert =p}
\E(\partial _{\alpha }f(\widehat{S}_{n}^{k}))\E(F_k^{\alpha })+\E\Big(U_{r}f\Big(\widehat{S}_{n}^{k},\frac{F_k}{n^{1/2}}\Big)\Big).
\end{equation*}%
We now use that $S_n^{k-1}=\widehat S_n^k+\frac{G_k}{n^{1/2}}$: the same reasoning for $G_{k}$ gives%
\begin{equation*}
\E(f(S_{n}^{k-1}))=\E(f(\widehat{S}_{n}^{k}))+\sum_{p=1}^{r}\frac{1}{p!n^{p/2}}%
\sum_{\left\vert \alpha \right\vert =p}\E(\partial _{\alpha }f(\widehat{S}%
_{n}^{k}))\E(G_{k}^{\alpha })+
\E\Big(U_{r}f\Big(\widehat{S}_{n}^{k},\frac{G_{k}}{n^{1/2}}\Big)\Big).
\end{equation*}%
By recalling that $\Delta_\alpha=\E(F^\alpha)-\E(G^\alpha)=0$ for $|\alpha|\leq 2$, the statement holds. $\square$

\medskip

Our aim is now to replace $\widehat{S}_{n}^{k}$\ by $S_{n}^{k-1}$ in the
development (\ref{taylor}). This opens the way to use a recurrence
procedure. 

\begin{lemma}\label{lemma-12}
Let $n\in \N,1\leq k\leq n$ and $r\in\N$. If $F\in L^{r+1}(\Omega)$ then for every  $f\in C^{r+2}_b(\R^{N})$ one has
$$
\E(f(S^{k}_n))-\E(f(S^{k-1}_n))
=\sum_{t=3}^r\frac 1{n^{t/2}}\,\E\big(\Psi_tf(S^{k-1}_n)\big)+\frac 1{n^{(r+1)/2}}\,\mathcal{R}_{r,n}^kf
$$
where $\Psi_t$ and $\mathcal{R}_{r,n}^k$ are defined in (\ref{Psi}) and (\ref{Rkf}), respectively.
\end{lemma}

\textbf{Proof.}
Consider the generical term $\E(\partial_\alpha f(\widehat S^{k}_n))$  of (\ref{taylor}). We recall that $\widehat S^{k}_n+G_k/\sqrt n=S^{k-1}_n$ and that $\widehat S^{k}_n$ and $G_k$ are independent. So, we apply (\ref{Taylor-back}) in Appendix \ref{BT} to $g(x)=\partial_\alpha f(\widehat S^{k}_n+x/\sqrt n)$ with $|\alpha|=p\leq r$, and we expand up to the maximum order $L$ such that $2L\leq r-p$. Hence we can write
\begin{align*}
\E(\partial_\alpha f(\widehat S^{k}_n))
&=\sum_{q=0}^{[(r-p)/2]}
\frac{(-1)^q}{2^q q! n^q}\sum_{|\beta|={2q}}
\theta_\beta\E(\partial_\beta\partial_\alpha f(S^{k-1}_n))
+\frac{1}{n^{[(r-p)/2]+1}}\widetilde U_{[\frac{r-p}2]}\partial_\alpha f\Big(S^{k-1}_n,\frac{G_k}{\sqrt n}\Big)
\end{align*}
where
$$
\widetilde
U_{L}g\Big(\widehat S^{k}_n,\frac{G_k}{\sqrt n}\Big)
=\frac{(-1)^{L+1}}{2^{L+1} L!}\sum_{|\beta|={2L+2}}\theta_\beta\int_0^1 s^L
\E\Big(\partial_\beta g\Big(\widehat S^{k}_n+\sqrt s\,\frac {G_k}{\sqrt n}\Big)\Big)ds
$$
By inserting in (\ref{taylor}) we get
\begin{align*}
\E(f(S^{k}_n))-\E(f(S^{k-1}_n))
=&\sum_{p=3}^r\frac 1{p! n^{p/2}}\sum_{|\alpha|=p}\Delta_\alpha
\sum_{q=0}^{[(r-p)/2]}
\frac{(-1)^q}{2^q q! n^q}\sum_{|\beta|={2q}}
\theta_\beta\E(\partial_\beta\partial_\alpha f(S^{k-1}_n))+\\
&+\sum_{p=3}^r\frac 1{p! n^{p/2}}\sum_{|\alpha|=p}\Delta_\alpha\frac{1}{n^{[\frac{r-p}2]+1}}\widetilde U_{[\frac{r-p}2]}\partial_\alpha f\Big(\widehat S^{k}_n,\frac{G_k}{\sqrt n}\Big)
+\frac 1{n^{(r+1)/2}}
\widetilde{\mathcal{R}}^k_{r,n}f\\
=&\sum_{p=0}^r\sum_{q=0}^{[(r-p)/2]}\frac  {(-1)^q}{2^qp!q! n^{(p+2q)/2}}\sum_{|\alpha|=p}
\sum_{|\beta|={2q}} \E\big(\partial_\beta\partial_\alpha f(S^{k-1}_n)\big)\theta_\beta\Delta_\alpha+\\
&+\frac 1{n^{(r+1)/2}}\mathcal{R}_{r,n}^kf
\end{align*}
in which, for the last line, we have used (\ref{Rkf}) -- recall that in the sum we can let the index $p$ start from 0 because as $p=0,1,2$, $\Delta_\alpha=0$. Now, by considering the change of variable $(t,s)=(p+2q,p)$ in the double sum above, we get
\begin{align*}
\E(f(S^{k}_n))-\E(f(S^{k-1}_n))
=&\sum_{t=0}^r\sum_{s=0}^ t
%
%
\frac  {(-1)^{\frac{t-s}2}}{2^{\frac{t-s}2}s!(\frac{t-s}2)! n^{t/2}}\sum_{|\alpha|=s}
\sum_{|\beta|={t-s}} \E\big(\partial_\beta\partial_\alpha f(S^{k-1}_n)\big)\theta_\beta\Delta_\alpha+\\
&+\frac 1{n^{(r+1)/2}}\mathcal{R}_{r,n}^kf\\
=&\sum_{t=0}^r\frac{1}{n^{t/2}}
\E\big(\Psi_tf(S^{k-1}_n)\big)+\frac 1{n^{(r+1)/2}}\mathcal{R}_{r,n}^kf.
\end{align*}
Since $\Psi_t=0$ for $t\leq 2$,  the statement holds. $\square$

\medskip

For $k=1,\ldots,n$, we define
\begin{equation}\label{Psik}
\mbox{$\Psi^{(1)}_t=\Psi_t$ and for $k\geq 2$, $\Psi^{(k)}_t=\Psi^{(k-1)}_t+\sum_{p=0}^t\Psi_p\Psi^{(k-1)}_{t-p}$, $t=0,1,\ldots$}
\end{equation}
Notice that $\Psi^{(k)}_t$ is a differential operator which is linked to the convolution w.r.t. $t$ between $\Psi_{\cdot}$ and the preceding operator $\Psi^{(k-1)}_{\cdot}$. We also notice that
$\Psi^{(k)}_t=0$ for $t=0,1,2$,
as an immediate consequence of the fact that $\Psi_t=0$ for $t\leq 2$. So, for $k\geq 2$ we can write
\begin{equation}\label{Psik-bis}
\Psi^{(k)}_t=\I_{\{t\geq 3\}}\Psi^{(k-1)}_t+\I_{\{t\geq 6\}}\sum_{p=3}^{t-3}\Psi_p\Psi^{(k-1)}_{t-p},\quad t=0,1,\ldots,
\end{equation}
We also define the following reminder operators: for $r\in\N$,
\begin{equation}\label{Phikr}
\Phi^{(k)}_{r,n}f=\sum_{j=1}^{k-1}\sum_{t=0}^r\mathcal{R}^{k-j}_{r-t,n}\Psi^{(j)}_tf+\mathcal{R}_{r,n}^kf.
\end{equation}
Note that, by definition, $\Phi^{(0)}_{r,n}=\mathcal{R}^0_{r,n}$ and $\Phi^{(k)}_{0,n}=\mathcal{R}^k_{0,n}$.

\begin{lemma}\label{lemma-bis}
Let $n\in \N,1\leq k\leq n$ and $r\in \N$. If $F\in L^{r+1}(\Omega)$ then for every  $f\in C^{r+2}_b(\R^{N})$ one has
$$
\E\big(f(S^{k}_n)\big)-
\E\big(f(S^{k-1}_n)\big)
=\sum_{t=3}^r\frac 1{n^{t/2}}\E\big(\Psi^{(k)}_tf(S^0_n)\big)+\frac 1{n^{(r+1)/2}}\Phi^{(k)}_{r,n}f,
$$
$\Psi^{(k)}_t$ and $\Phi^{(k)}_{r,n}$ being given in (\ref{Psik}) and (\ref{Phikr}), respectively.

\end{lemma}

\proof
We consider the development in Lemma \ref{lemma-12}:
$$
\E(f(S^{k}_n))-\E(f(S^{k-1}_n))
=\sum_{t=0}^r\frac 1{n^{t/2}}\,\E\big(\Psi_tf(S^{k-1}_n)\big)+\frac 1{n^{(r+1)/2}}\,\mathcal{R}_{r,n}^kf.
$$
For $t\leq r$, we apply such development up to order $r-t$ to $\E\big(\Psi_tf(S^{k-1}_n)\big)$ and we get
$$
\E(\Psi_tf (S^{k-1}_n))=\E(\Psi_tf(S^{k-2}_n))
+\sum_{p=0}^{r-t}\frac 1{n^{p/2}}\,\E\big(\Psi_p\Psi_tf(S^{k-2}_n)\big)+\frac 1{n^{(r-t+1)/2}}\,\mathcal{R}^{k-1}_{r-t,n}\Psi_tf.
$$
By inserting, we obtain
\begin{align*}
\E(f(S^{k}_n))-\E(f(S^{k-1}_n))
=&\sum_{t=0}^r\frac 1{n^{t/2}}\E(\Psi_tf(S^{k-2}_n))
+\sum_{t=0}^r\sum_{p=0}^{r-t}\frac 1{n^{(t+p)/2}}\,\E\big(\Psi_p\Psi_tf(S^{k-2}_n)\big)+\\
&+\frac 1{n^{(r+1)/2}}\,
\sum_{t=0}^r
\mathcal{R}^{k-1}_{r-t,n}\Psi_tf+\frac 1{n^{(r+1)/2}}\,\mathcal{R}_{r,n}^kf
\end{align*}
and by a change of variable in the second sum above we get
\begin{align*}
\E(f(S^{k}_n))-\E(f(S^{k-1}_n))
=&\sum_{t=0}^r\frac 1{n^{t/2}}\E(\Psi^{(2)}_tf(S^{k-2}_n))
+\frac 1{n^{(r+1)/2}}\,\Big[\sum_{t=0}^r
\mathcal{R}^{k-1}_{r-t,n}\Psi_tf+\mathcal{R}_{r,n}^kf\Big].
\end{align*}
By iterating the same procedure up to step $k$, we obtain the statement. $\square$

\medskip

\def\mU{\mathcal U}
We now set
\begin{equation}\label{TR}
T^n_t=\sum_{k=1}^n\Psi^{(k)}_t\quad\mbox{and}\quad
\mU ^n_r=\sum_{k=1}^n\Phi^{(k)}_{r,n}
\end{equation}
$\Psi_t^{(k)}$ and $\Phi^{(k)}_{r,n}$ being given in (\ref{Psik}) and (\ref{Phikr}), respectively.

\begin{proposition}\label{prop1}
Let $n\in \N,1\leq k\leq n$ and $r\in\N$. If $F\in L^{r+1}(\Omega)$ then for every  $f\in C^{r+2}_b(\R^{N})$ one has
$$
\E\big(f(S^{n}_n)\big)-
\E\big(f(S^{0}_n)\big)
=\sum_{t=3}^r\frac 1{n^{t/2}}\E\big(T_t^nf(S^0_n)\big)+
\frac 1{n^{(r+1)/2}}\mU^n_r f,
$$
where
$T^n_t$ and $\mU^n_r$ are defined in (\ref{TR}).
\end{proposition}

\proof Since $\E\big(f(S^{n}_n)\big)-\E\big(f(S^{0}_n)\big)=\sum_{k=1}^n\big(\E\big(f(S^{k}_n)\big)-
\E\big(f(S^{k-1}_n)\big)
\big)$, the statement immediately follows from Lemma \ref{lemma-bis}. $\square$

\medskip

We give now an explicit expression for the operators $\Psi^{(k)}_t$ in (\ref{Psik}) and, as a consequence, for $T^n_t$ in (\ref{TR}).
For $\Psi_t$ given in (\ref{Psi}), $i=1,2,\ldots$, we set
$$
\mA^1_t=\Psi_t\quad\mbox{and for}\quad i\geq 1, \quad \mA^{i+1}_{t}
=\sum_{p=0}^{t}
\Psi_{p}\mA^i_{t-p}.
$$
Since $\Psi_t=0$ for $t=0,1,2$, straightforward computations give that
$\mA^i_t=0$ if $t<3i$, so that we can also write
\begin{equation}\label{Ai-rec}
\mA^1_t=\Psi_t\quad\mbox{and for}\quad i\geq 1, \quad \mA^{i+1}_{t}
=\sum_{p=3}^{t-3i}
\Psi_{p}\mA^i_{t-p}.
\end{equation}
We can give an alternative representation for the $\mA^i_t$'s.
We set $\mM$ the set of all multiindexes and for $\alpha,\beta\in\mM$ (possibly with different length), we set $(\alpha,\beta)\in\mM$ the associated concatenation. So, for $\gamma\in\mM$ we define
$$
A_\gamma=\{(\alpha,\beta)\,:\,(\alpha,\beta)=\gamma\}
$$
and
\begin{equation}\label{ci}
c^1_\gamma
=\sum_{(\alpha,\beta)\in A_\gamma}
\frac{(-1)^{\frac{|\beta|}2}}{2^{\frac{|\beta|}2}|\alpha|!(\frac{|\beta|}2)!}\Delta_\alpha\theta_\beta \quad\mbox{and for $i\geq 1$,}\quad
c^{i+1}_\gamma=\sum_{(\alpha,\beta)\in A_\gamma}c^1_\alpha c^i_\beta,\quad i\geq 1.
\end{equation}

%
Since $c^1_\gamma=0$ if $|\gamma|<3$, by recurrence one gets $c^i_\gamma=0$ if $|\gamma|<3i$ for every $i$. Then, straightforward computations give that, for $i\geq 1$,
\begin{equation}\label{Aitbis}
\mA^i_t=\sum_{\gamma\,:\,|\gamma|=t}c^i_\gamma\partial_\gamma,\quad\mbox{with $\{c^i_\gamma\}_{\gamma\in\mM}$ given in (\ref{ci})}.
\end{equation}
It is immediate to see that for every $\gamma\in\mM$ there exists $C$ such that for every $i\geq 1$
\begin{equation}\label{ci-est}
|c^i_\gamma|\leq C \sup_{|\alpha|\leq |\gamma|}|\Delta_\alpha|^i.
\end{equation}
As a consequence, for  $t, q\geq 0$ there exists $C>0$ (depending on $t, q$ only) such that
for every $i\geq 1$  and $f\in C^{t+q}_b(\R^N)$
\begin{equation}\label{Ait-est}
\|\mA^{i}_{t}f\|_{q,\infty}
\leq C\sup_{|\alpha|\leq t}|\Delta_\alpha|^i\times \|f\|_{t+q,\infty}
\leq C(1+\E(|F|^t))^{i-1}\sup_{|\alpha|\leq t}|\Delta_\alpha|\times \|f\|_{t+q,\infty}.
\end{equation}
Moreover, the $\mA^i_t$'s give the following representation formula for the $\Psi^{(k)}_t$'s:
\begin{proposition}\label{miracle-real}
For every $k\geq 1$ the operator $\Psi^{(k)}$  given in (\ref{Psik}) can be written as
$$
\Psi^{(k)}_t=\sum_{i=1}^{[t/3]}
Q_{i-1}(k)\mA^i_t,\quad t=0,1,\ldots
$$
where $Q_{i-1}(k)$  is defined as follows:
$$
\mbox{$Q_0(k)=1$ and for $l\geq 1$, $Q_l(k)=\displaystyle\sum_{j=l+1}^kQ_{l-1}(j-1)$.}
$$
In particular, $Q_l(k)=0$ if $k\leq l$ and $Q_l(k)>0$ otherwise.
\end{proposition}

\proof We have already observed that if $[t/3]=0$ then $\Psi^{(k)}_t=\Psi_t=0$ for every $k$ and if $[t/3]=1$ then $\Psi^{(k)}_t=\Psi_t$ for every $k$, see (\ref{Psik-bis}), so the formulas agree. We now assume that the formula is true for $[t/3]=j\geq 1$ and for every $k$, and we prove it for $[t/3]=j+1$ and for every $k$. We recall that $\Psi^{(k)}_t=\Psi^{(k-1)}_t+\sum_{p=3}^{t-3}\Psi_p\Psi^{(k-1)}_{t-p}$. But if $[t/3]=j+1$ then $[(t-p)/3]\leq j$ for any $p=3,\ldots,t-3$, so that by induction $\Psi^{(k-1)}_{t-p}$ fulfils the formula. Therefore, we can write
\begin{align*}
\Psi^{(k)}_t
=&\Psi^{(k-1)}_t
+\sum_{p=3}^{t-3}
\sum_{i=1}^{[(t-p)/3]}
Q_{i-1}(k-1)
\Psi_{p}A^i_{t-p}.
\end{align*}
We do a change of variable in the last sum:
the condition $i\leq [(t-p)]/3$ gives $3i\leq t-p$, that is $p\leq t-3i$, and if $p\geq 3$ then $i\leq [t/3]-1$. So, by using also (\ref{Ai-rec}) we get
\begin{align*}
\Psi^{(k)}_t
-\Psi^{(k-1)}_t
=&\sum_{i=1}^{[t/3]-1}
Q_{i-1}(k-1)
\sum_{p=3}^{t-3i}
\Psi_{p}\mA^i_{t-p}
=
\sum_{i=1}^{[t/3]-1}
Q_{i-1}(k-1)
\mA^{i+1}_{t}
=\sum_{i=2}^{[t/3]}
Q_{i-2}(k-1)\mA^{i}_{t}.
\end{align*}
By summing
\begin{align*}
\Psi^{(k)}_t
=&\Psi_t
+\sum_{i=2}^{[t/3]}
\sum_{j=2}^kQ_{i-2}(j-1)
\mA^i_t
=Q_0(k)\mA^1_t +\sum_{i=2}^{[t/3]}
\sum_{j=2}^kQ_{i-2}(j-1)
\mA^i_t
\end{align*}
and the statement holds for $Q_0(k)=1$ and $Q_{i-1}(k)=\sum_{j=2}^kQ_{i-2}(j-1)$, $i\geq 2$. We now prove that $Q_l(k)=0$ if $k\leq l$ and
$Q_l(k)>0$ for $k\geq l+1$. For $l=1$, $Q_l(k)=k-1$, and the statement holds. If we assume that $Q_l(k)$ is not null for $k\geq l+1$ then
$$
Q_{l+1}(k)=\sum_{j=2}^kQ_{l}(j-1)\I_{j-1\geq l+1}=\sum_{j=2}^kQ_{l}(j-1)\I_{j\geq l+2}
$$
and this is null for $k\leq l+1$ and  strictly positive if $k\geq l+2$.
 $\square$

\medskip

We now give an explicit formula for $T^n_t$, namely we write it in such a way that $n\mapsto T^n_t$ is a polynomial whose coefficients will be explicitly written.
To this purpose, we need to handle polynomials of the type
$$
n\mapsto \SS_l(n-1)
=\sum_{k=1}^{n-1}k^l, \quad l\in\N, n\geq 1.
$$
We recall the exact expansion for $\SS_l(L)=\sum_{k=1}^{L}k^l$:
\begin{equation}\label{SL}
\SS_l(L)
=\frac 1{l+1}\sum_{p=1}^{l+1}\coeffbin{l+1}{p}B_{l+1-p}\,L^p
\end{equation}
where $\{B_m\}_m$ denotes the sequence of the (second) Bernoulli numbers (which are in fact defined as the numbers for which the above equality holds, see \cite{[AS]}), whose first numbers are given by
$$
B_0=1,\  B_1=\frac 12,\  B_2=\frac 16,\  B_3=0,\  B_4=-\frac 1{30},\  B_5=0,\  B_6=\frac 1{42},\  B_7=0,\  B_8=-\frac 1{30},\ldots
$$
Then, straightforward computations give that  for $l\in\N$ and $n\geq 1$,
$$
\SS_l(n-1)=\sum_{k=1}^{n-1}k^l=\sum_{q=0}^{l+1}b_{l, q}n^q
$$
where the sequence $(b_{l,q})_{q=0,\ldots,l+1}$ is given by
\begin{equation}\label{Blq}
b_{l,q}=
\frac 1{l+1}\sum_{p=q\vee 1}^{l+1}\coeffbin{l+1}{p}B_{l+1-p}\coeffbin{p}{q}(-1)^{p-q},\quad
q=0,1,\ldots,l+1\mbox{ and } l\in\N,
\end{equation}
in which $B_l$, ${l\geq 0}$, denote the (second) Bernoulli numbers.
Just as an example:
$$
\hskip -3.cm
\begin{array}{lllll}
\bullet \ l=0:\qquad & b_{0,0}=-1, & b_{0,1}=1; & &\smallskip\\
\bullet \ l=1: & b_{1,0}=0, &b_{1,1}=-\frac 12, & b_{1,2}=\frac 12; &\smallskip\\
\bullet \ l=2: & b_{2,0}=0, &b_{2,1}=\frac 16, &b_{2,2}=-\frac 12, & b_{2,3}=\frac 13.
\end{array}
$$

Then one has the following.
\begin{proposition}\label{prop-Tn}
Let $n\geq 1$, $r\in \N$ and $F\in L^{q_r+1}(\Omega)$, where $q_r=\max(r,2)$. For $t\leq r$, let $T^n_t$ be defined as in (\ref{TR}). Then,
$$
T^{n}_t=\sum_{i=1}^{[t/3]}
P_{i}(n)\mA^i_t,\quad t=0,1,\ldots
$$
where $P_i(n)=0$ if $n<i$ and for $n\geq i$,
\begin{equation}\label{Pin}
P_i(n)=\sum_{p=0}^i a_{i,p}n^p,\quad i= 1,\ldots,n
\end{equation}
with
\begin{equation}\label{aip}
\begin{array}{l}
a_{1,0}=0,\quad a_{1,1}=1\quad\mbox{and for $i\geq 1$}\smallskip\\
\displaystyle
a_{i+1,0}
=\sum_{l=0}^i
a_{i,l}b_{l, 0}
-\sum_{l=0}^i
a_{i,l}\SS_l(i-1),
\quad
a_{i+1,p}
=\sum_{l=p-1}^i
a_{i,l}b_{l, p},\quad p=1,\ldots,i
\end{array}
\end{equation}
the sequence $(b_{l,p})_{p=0,\ldots,l+1}$ being defined in (\ref{Blq}) and $\SS_l(i-1)$ being given in (\ref{SL}).
\end{proposition}
\textbf{Proof.}
Since $T^n_t=\sum_{k=1}^n\Psi^{(k)}_t$, we get
$$
T^n_t=\sum_{i=1}^{[t/3]}\sum_{k=1}^nQ_{i-1}(k)\mA^i_t
$$
so that $P_i(n)=\sum_{k=1}^nQ_{i-1}(k)=\sum_{j=2}^{n+1}Q_{i-1}(j-1)=Q_i(n+1)$. As a consequence, $P_i(n)=0$ if $n+1\leq i$, that is $n<i$. So, let $n\geq i$.
We have $P_1(n)=\sum_{k=1}^nQ_0(k)=n$ and for $i\geq 2$,
\begin{equation}\label{rec-Pi}
P_i(n)
=Q_i(n+1)
=\sum_{j=2}^{n+1}Q_{i-1}(j-1)\I_{j-1\geq i}
=\sum_{k=i-1}^{n-1}Q_{i-1}(k+1)
=\sum_{k=i-1}^{n-1}P_{i-1}(k).
\end{equation}
Since $P_1(n)=n$, we get $a_{1,0}=0$ and $a_{1,1}=1$. In order to compute the sequence
$(a_{i,l})_{l=0,\ldots,i}$, we use a recurrence argument. For $i\geq 1$, one has
\begin{align*}
P_{i+1}(n)=
&\sum_{k=i}^{n-1}P_{i}(k)
=\sum_{k=i}^{n-1}
\sum_{l=0}^i a_{i,l}k^l
=\sum_{l=0}^i
a_{i,l}
\sum_{k=i}^{n-1}k^l
=\sum_{l=0}^i
a_{i,l}\big(\SS_l(n-1)-\SS_l(i-1)\big)\\
=
&\sum_{l=0}^i
a_{i,l}\SS_l(n-1)
-\sum_{l=0}^i
a_{i,l}\SS_l(i-1)
=\sum_{l=0}^i
a_{i,l}\sum_{p=0}^{l+1}b_{l, p}n^p
-\sum_{l=0}^i
a_{i,l}\SS_l(i-1)\\
=&
\sum_{p=0}^{i+1}n^p
\sum_{l=0\vee(p-1)}^i
a_{i,l}b_{l, p}
-\sum_{l=0}^i
a_{i,l}\SS_l(i-1)
\end{align*}
and (\ref{aip}) follows.
$\square$

\medskip

We are now ready to prove our result on the asymptotic expansion for smooth functions.
We set:

\smallskip

$\bullet$ for $m\geq 1$ and $f\in C^{m}_b(\R^N)$,
\begin{equation}\label{D}
\mD_m f
=\sum_{\mbox{\scriptsize{$\begin{array}{c}t=3\vee m\\t-m \mbox{ even}\end{array}$}}}^{3m}
\sum_{i=1\vee \frac{t-m}2}^{[t/3]}a_{i,\frac{t-m}2}
\E\big(\mA^i_t f(G)\big);
\end{equation}
$\bullet$ for $r\geq 2$ and $f\in C^{r+2}_b(\R^N)$,
\begin{equation}\label{mE}
\mE^n_{r}f
=n^{\frac{[r/3]+1}2}\times\Big[\sum_{m=[\frac r3]+1}^{r}\frac 1{n^{\frac m2}}\sum_{\mbox{\scriptsize{$\begin{array}{c}t=3\vee m\\t-m \mbox{ even}\end{array}$}}}^{(3m)\wedge r}
\sum_{i=1\vee \frac{t-m}2}^{[t/3]}a_{i,\frac{t-m}2}
\E\big(\mA^i_t f(G)\big)
+\frac 1{n^{\frac{r+1}2}}\mU^n_r f\Big].
\end{equation}

Then we have

\begin{theorem}\label{dev1}
Let $r\geq 2$. If $F\in L^{r+1}(\Omega)$,  then for every $f\in C^{r+3}_b(\R^N)$ one has
\begin{align*}
\E\big(f(S_n)\big)-
\E\big(f(G)\big)
=&\sum_{m=1}^{[r/3]}\frac 1{n^{\frac m2}}\mD_m f+
\frac 1{n^{\frac{[r/3]+1}2}}\mE^n_{r}f
\end{align*}
where $\mD_m f$ and $\mE^n_rf$ are defined in (\ref{D}) and (\ref{mE}), respectively. \end{theorem}

\begin{remark}
At this stage, we could prove that
\begin{equation}\label{est-nohere}
|\mathcal{E}^n_{r}f|
\leq C(1+\E(|F|^{r+1}))^{[r/3]\vee 1}\Big[\|f\|_{r+3,\infty}\sup_{|\alpha|\leq r}|\Delta_\alpha|
+\|f\|_{r+2,\infty}\frac 1{n^{\frac{r-[r/3]-2}2}}\Big],
\end{equation}
$C$ denoting a suitable constant depending on $r$ and $N$ only. But since we aim to deal with the distance in total variation, we need a representation and an estimate of the reminder in terms of $f$ and not of its derivatives. So, we skip this point and we postpone the problem to next section.
\end{remark}

\textbf{Proof of Theorem \ref{dev1}.}
Take $r\geq 2$. We use Proposition \ref{prop1}: for every $n\in \N$ and $%
f\in C^{r+2}_b(\R^{N})$ we have
\begin{align*}
\E\big(f(S^{n}_n)\big)-
\E\big(f(S^{0}_n)\big)
=&\sum_{t=3}^r\frac 1{n^{\frac t2}}\sum_{i=1}^{[t/3]}
P_{i}(n)\E\big(\mA^i_t f(G)\big)+
\frac 1{n^{(r+1)/2}}\mU^n_r f\\
=&\sum_{t=3}^r\frac 1{n^{\frac t2}}\sum_{i=1}^{[t/3]}
\sum_{p=0}^i a_{i,p}n^p\,
\E\big(\mA^i_t f(G)\big)+
\frac 1{n^{(r+1)/2}}\mU^n_r f\\
=&\sum_{t=3}^r\sum_{p=0}^{[t/3]}\frac 1{n^{\frac t2-p}} \sum_{i=1\vee p}^{[t/3]}
a_{i,p}
\E\big(\mA^i_t f(G)\big)+
\frac 1{n^{(r+1)/2}}\mU^n_r f.
\end{align*}
So, by recalling that $S_n=S_n^n$ and $G\stackrel{\mathcal{L}}{=}S_n^0$
we obtain
\begin{align*}
\E\big(f(S_n)\big)-
\E\big(f(G)\big)
=&\sum_{t=3}^r\sum_{p=0}^{[t/3]}\frac 1{n^{\frac t2-p}} \sum_{i=1\vee p}^{[t/3]}a_{i,p}
\E\big(\mA^i_t f(G)\big)+
\frac 1{n^{(r+1)/2}}\mU^n_r f
\end{align*}
We set now $t-2p=m$, so $t-m$ is an even number. Now, $p\geq 0$ gives that $t\geq m$ and since $t\geq 3$ then $t\geq 3\vee m$ and $m\leq r$; $p\leq [t/3]$ gives that $(t-m)/2\leq [t/3]$. Therefore, the sum over $t\leq r$ must be done on the set $\{t\,:\, 3\vee m\leq t\leq r, t-m\mbox{ even, } t-2[t/3]\leq m\}$. It is easy to see that this set equals to
$\{t\,:\, 3\vee m\leq t\leq (3m)\wedge r, t-m\mbox{ even}\}$.
So, we obtain
\begin{align*}
\E\big(f(S_n)\big)-
\E\big(f(G)\big)
=&\sum_{m=1}^r\frac 1{n^{\frac m2}}\sum_{\mbox{\scriptsize{$\begin{array}{c}t=3\vee m\\t-m \mbox{ even}\end{array}$}}}^{(3m)\wedge r}
\sum_{i=1\vee \frac{t-m}2}^{[t/3]}a_{i,\frac{t-m}2}
\E\big(\mA^i_t f(G)\big)+
\frac 1{n^{(r+1)/2}}\mU^n_r f.
\end{align*}
The statement now follows by using (\ref{D}) (notice that $3m\leq r$ if $m\leq [r/3]$) and (\ref{mE}).
$\square$\medskip

\subsection{Regularized functions and estimate of the reminder}\label{sect-regularization}

Our problem is now to prove an estimate for the reminder in the development for a function $f$ in terms of $\left\Vert f\right\Vert _{\infty }$ instead of $\left\Vert f\right\Vert _{r+1,\infty }$. To this purpose, we need some preliminary results.

For $\delta >0$, we denote by $\gamma _{\delta }$ the
density of the centred Gaussian law in $\R^N$ of variance $\delta I$ and for $%
f:\R^{N}\rightarrow \R$ we denote $f_{\delta }=f\ast \gamma _{\delta }.$ Using
standard integration by parts on $\R^{N}$, one may prove that for each $r\in \N$
there exists an universal constant $C$ (depending on $N$ and $r$ only) such
that for every multiindex $\alpha $ with $\left\vert \alpha \right\vert =r$
one has
\begin{equation}
\left\Vert \partial _{\alpha }f_{\delta }\right\Vert _{\infty }\leq \frac{C}{%
\delta ^{r/2}}\left\Vert f\right\Vert _{\infty }.  \label{abs26}
\end{equation}%
We give now some estimates following from Lemma \ref{Lemma2.5} with $\Theta=1$, which is actually Lemma 2.5 in \cite{[BC]}.

\begin{lemma}\label{pippo}
Suppose that $\mu _{F}\succeq \leb_{N}.$
There exist universal constants $C>0$ and $b>4$, depending on $N$ only, such that for every $\delta >0$, $n\in\N$ and
for every bounded and measurable function $f:{\mathbb{R}}^{N}\rightarrow {\mathbb{R}}$ one has
\begin{equation}\label{abs30'}
\left\vert {\mathbb{E}}(f(S_n))-{\mathbb{E}}(f_{\delta
}(S_n))\right\vert
\leq C\left\Vert f\right\Vert _{\infty }(1+\E(|F|))\Big(e^{-n/C}+\delta ^{1/b} n^{(b-2)/(2b)}\Big).
\end{equation}
\end{lemma}

\textbf{Proof.} Let $K\geq 1$ and $\Psi _{K}\in C^{\infty }(\R^{N})$ be such
that $1_{B_{K}(0)}\leq \Psi _{K}\leq 1_{B_{K+1}(0)}$ and such that, for some $L>0$,
$\left\Vert \partial _{\alpha }\Psi _{K}\right\Vert _{\infty }\leq L$ for
every multiindex $\alpha $. Then we have
\[
\left\vert {\mathbb{E}}(f(S_n))-{\mathbb{E}}(f(\Psi _{K}(S_n)S_n))\right\vert
\leq
\left\Vert f\right\Vert _{\infty }\P(\left\vert S_n\right\vert \geq K)
\leq
\left\Vert f\right\Vert _{\infty }\frac{\E(|S_n|)}{K}
\leq
\left\Vert f\right\Vert _{\infty }\frac{\sqrt n}{K}\, \E(|F|)
\]%
and in a similar way $\left\vert {\mathbb{E}}(f_{\delta }(S_n))-{%
\mathbb{E}}(f_{\delta }(\Psi _{K}(S_n)S_n))\right\vert\leq \left\Vert f\right\Vert
\E(|F|)\sqrt n/K$. So we can write
\begin{align*}
|\E(f(S_n))-\E(f_\delta(S_n))|
\leq & |\E(f(S_n))-\E(f(\Psi _{K}(S_n)S_n))|
+|\E(f_\delta(S_n))-\E(f_\delta(\Psi _{K}(S_n)S_n))|\\
&+|\E(f(\Psi _{K}(S_n)S_n))-\E(f_\delta(\Psi _{K}(S_n)S_n))|\\
\leq &
2\E(|F|)\left\Vert f\right\Vert _{\infty }\,\frac{\sqrt n}{K}
+|\E(f(\Psi _{K}(S_n)S_n))-\E(f_\delta(\Psi _{K}(S_n)S_n))|.
\end{align*}
As for the last term in the above right hand side, we apply Lemma \ref{Lemma2.5} with $\Theta=1$  and $X=\Psi_K(S_n)S_n$: there exist some universal constants $C,p,a$ depending only on $N$ such that for every $\varepsilon >0,\delta >0$ and every $f\in L^\infty(\R^N)$ then
\begin{align*}
&\left\vert {\mathbb{E}}(f(\Psi_K(S_n)S_n))-{\mathbb{E}}(f_{\delta
}(\Psi_K(S_n)S_n))\right\vert
\leq C\left\Vert f\right\Vert _{\infty }\times\\
&\qquad\qquad\times\Big({\mathbb{P}}%
(\det\sigma _{\Psi_K(S_n)S_n}<\varepsilon )+\frac{\sqrt{\delta }}{\varepsilon ^{p}}%
(1+\left\Vert \Psi_K(S_n)S_n\right\Vert _{3,p}+\left\Vert L(\Psi_K(S_n)S_n)\right\Vert _{1,p})^{a}\Big).
\end{align*}
We note that we are forced to introduce the localization $\Psi_K(S_n)$ because in the above estimate it appears $\|\Psi_K(S_n)S_n\|_p$ with $p>1$: since the r.v.'s are only square integrable, if we take $\Psi_K\equiv 1$ then in principle we do not know if such norm is finite.

Now, on the set $\{|S_n|\leq K\}$ we have $\det\sigma _{\Psi_K(S_n)S_n}=\det\sigma _{S_n}$, so that
\begin{align*}
{\mathbb{P}}%
(\det\sigma _{\Psi_K(S_n)S_n}<\varepsilon )
&\leq {\mathbb{P}}%
(\det\sigma _{S_n}<\varepsilon)+
{\mathbb{P}}%
(|S_n|>K )
\leq {\mathbb{P}}%
(\det\sigma _{S_n}<\varepsilon)+\frac{\E(|S_n|)}K\\
&\leq {\mathbb{P}}%
(\det\sigma _{S_n}<\varepsilon)+\E(|F|)\frac{\sqrt n}K.
\end{align*}
By taking $\varepsilon= \varepsilon_*/2$ as in Lemma \ref{est-sigma}, (\ref{abs29}) gives
\begin{align*}
{\mathbb{P}}%
(\det\sigma _{\Psi_K(S_n)S_n}<\varepsilon )
&\leq Ce^{-n/C}+\E(|F|)\frac{\sqrt n}K.
\end{align*}
Therefore, we can write
\begin{align*}
&\left\vert {\mathbb{E}}(f(S_n))-{\mathbb{E}}(f_{\delta
}(S_n))\right\vert
\leq C\left\Vert f\right\Vert _{\infty }\times\\
&\qquad\qquad\times\Big(e^{-n/C}+\E(|F|)\,\frac{\sqrt n}K+\sqrt \delta \big(1+\left\Vert \Psi_K(S_n)S_n\right\Vert _{3,p}+\left\Vert L(\Psi_K(S_n)S_n)\right\Vert _{1,p}\big)^{a}\Big).
\end{align*}
We use now Lemma \ref{estK} in Appendix \ref{app-sob}: inequalities (\ref{estK-1}) and (\ref{estK-2}) give
\begin{align*}
\left\Vert \Psi_K(S_n)S_n\right\Vert _{3,p}+\left\Vert L(\Psi_K(S_n)S_n)\right\Vert _{1,p}
&\leq
CK\big(1+\|S_n\|_{1,3,4p}\big)^{6}
+
CK\big(1+\|S_n\|_{1,2,8p}\big)^{5}\big(1+
\|LS_n\|_{1,4p}\big)\\
&\leq
CK\big(1+\|S_n\|_{1,3,8p}+\|LS_n\|_{1,4p}\big)^{6}.
\end{align*}
By using (\ref{abs18}) and (\ref{abs19}) we have
\begin{align*}
\left\Vert \Psi_K(S_n)S_n\right\Vert _{3,p}+\left\Vert L(\Psi_K(S_n)S_n)\right\Vert _{1,p}
&\leq CK,
\end{align*}
so that
\begin{align*}
\left\vert {\mathbb{E}}(f(S_n))-{\mathbb{E}}(f_{\delta
}(S_n))\right\vert
&\leq C\left\Vert f\right\Vert _{\infty }\Big(e^{-n/C}+\E(|F|)\,\frac{\sqrt n}K+\sqrt{\delta }\, K^{a}\,\Big)\\
&\leq C\left\Vert f\right\Vert _{\infty }(1+\E(|F|))\Big(e^{-n/C}+\frac{\sqrt n}K+\sqrt{\delta }\, K^{a}\,\Big).
\end{align*}
We now optimize on $K$ by taking it in order that $\sqrt n/K=\sqrt{\delta }\, K^a$.
Straightforward computations give now (\ref{abs30'}), with $\frac 1b=\frac 12(1-\frac a{a+1})<\frac 14$.
$\square$

\begin{remark}\label{noId2}
We stress that when $C(F)\neq Id$ then the constant in (\ref{abs29}) depends on $\underline{\lambda}(F)$. As a consequence, this dependence holds for the constant $C$ appearing in (\ref{abs30'}) as well.

\end{remark}

We now propose the following key result, allowing us to deal with the remaining terms.

\begin{lemma}\label{4.4-new}
Suppose that $\mu _{F}\succeq \leb_{N}$.  Let $\alpha$ and $\beta$ denote multiindexes, with $|\alpha|=r$ and $|\beta|=m$. If $F\in L^m(\Omega)$, then there exists a constant $C$ (which depends on $N$, $r$ and $m$) such that for every $f\in L^{\infty}(\R^N)$, $\delta>0$, $n\geq 1$ and $\lambda\in\R$ then
\begin{align*}
&\Big|\E\Big(\partial _{\alpha }f_{\delta }\Big(\widehat{S}%
_{n}^{k}+\lambda \frac{F_k}{n^{1/2}}\Big)F_k^{\beta }\Big)\Big|
\leq C\|f\| _{\infty }\E(|F|^m)\big(1+\delta ^{-r/2} e^{-n/C}\big),\\
&\Big|\E\Big(\partial _{\alpha }f_{\delta }\Big(\widehat{S}%
_{n}^{k}+\lambda \frac{G_{k}}{n^{1/2}}\Big)G_k^{\beta }\Big)\Big|
\leq C\|f\| _{\infty }\E(|G|^m)\big(1+\delta ^{-r/2} e^{-n/C}\big),
\end{align*}
in which $f_{\delta }=f\ast \gamma _{\delta }$, $\gamma_\delta$ being the centred normal density in $\R^N$ with covariance matrix $\delta I$.

\end{lemma}

\textbf{Proof.} Without loss of generality, we suppose that $n$ is even and we study separately the cases $k\leq n/2$ and $k\geq n/2+1$ - if $n$ was odd, it would be sufficient to study $k\leq (n-1)/2$ and $k\geq (n-1)/2+1$.

\smallskip

\textbf{Case 1:} $k\leq n/2.$ We denote
\begin{equation*}
A_k=\frac{1}{n^{1/2}}\Big(\sum_{i=1}^{k-1}F_i+\sum_{i=k+1}^{n/2}G_{i}\Big)+\lambda \frac{F_k}{n^{1/2}},\qquad B=\frac{1}{n^{1/2}}\sum_{i=n/2+1}^{n}G_{i}
\end{equation*}%
so that
$$
\widehat{S}_{n}^{k}+\lambda\frac{F_k}{n^{1/2}}=A_k+B.
$$
Notice that $B$ is a Gaussian random variable with covariance $\frac{1}{2}I$
which is independent of $A_k$ and of $F_k.$ Using integration by parts
with respect to $B$ we may find a random variable $H_{\alpha}$ having all moments and
\begin{equation*}
\E\Big(\partial _{\alpha }f_{\delta }\Big(\widehat{S}%
_{n}^{k}+\lambda \frac{F_k}{n^{1/2}}\Big)F_k^{\beta }\Big)
=\E\big(\partial
_{\alpha }f_{\delta }(A_k+B)F_k^{\beta }\big)=\E(f_{\delta
}(A_k+B)F_k^{\beta }H_{\alpha}).
\end{equation*}%
Since $F_k$ and $H_{\alpha}$ are independent, $H_\alpha$ being a suitable function of $G_{n/2},\ldots,G_n$, it follows that
\begin{equation*}
\Big| \E\Big(\partial _{\alpha }f_{\delta }\Big(\widehat{S}%
_{n}^{k}+\lambda\frac{F_k}{n^{1/2}})F_k^{\beta }\Big)
\Big| \leq
C\left\Vert f_{\delta }\right\Vert _{\infty }\E(|F_k|^m)
\E(|H_\alpha|)
\leq C\left\Vert f\right\Vert _{\infty
}\E(|F|^m)
\end{equation*}%
Similarly, we obtain
\begin{equation*}
\Big| \E\Big(\partial _{\alpha }f_{\delta }\Big(\widehat{S}%
_{n}^{k}+\lambda\frac{G_{k}}{n^{1/2}})G_k^{\beta }\Big)
\Big| \leq
C\left\Vert f\right\Vert _{\infty
}\E(|G|^m).
\end{equation*}%

\textbf{Case 2:} $k>n/2.$ We denote
\begin{equation*}
A=\frac{1}{n^{1/2}}\sum_{i=1}^{n/2}F_i,\qquad
B_k=\frac{1}{n^{1/2}}%
(\sum_{i=n/2+1}^{k-1}F_i+\sum_{i=k+1}^{n}G_{i})+\lambda\frac{F_k}{n^{1/2}}
\end{equation*}%
so that
$$
\widehat{S}_{n}^{k}+\lambda\frac{F_k}{n^{1/2}}=A+B_k.
$$
We notice that
$$
A=\frac 1{\sqrt{2}} S_{n/2},
$$
so we can use the noise from the absolutely continuous r.v.'s $V_1,\ldots,V_{n/2}$ ``inside'' $S_{n/2}$, as already seen in Section \ref{AMC}. 
We then proceed to use integration by parts w.r.t. the noise from $A$.

We notice that $\sigma_A=\frac 12\sigma_{S_{n/2}}$ and that the covariance matrix $\sigma_{S_{n/2}}$ of $S_{n/2}$ may degenerate. So, we use a localization: we consider a function $\phi \in C^{1}(\R_{+})$ such that $%
\I_{(\varepsilon _{\ast }/2,\infty )}\leq \phi \leq \I_{(\varepsilon _{\ast
},\infty )}$ and $\left\Vert \nabla \phi \right\Vert _{\infty }\leq
2/\varepsilon _{\ast }$ with $\varepsilon _{\ast }$ given in (\ref{abs28})$.$
Then we write%
\begin{equation*}
\E\Big(\partial _{\alpha }f_{\delta }\Big(\widehat{S}%
_{n}^{k}+\lambda \frac{F_k}{n^{1/2}}\Big)F_k^{\beta }\Big)
=\E\big(\partial _{\alpha }f_{\delta }(A+B_k)F_k^{\beta }\big)=I+J
\end{equation*}%
with%
\begin{eqnarray*}
I &=&\E\big(\partial _{\alpha }f_{\delta
}(A+B_{k})F_k^{\beta }\phi (\det \sigma _{A})\big), \\
J &=&\E\big(\partial _{\alpha }f_{\delta
}(A+B_{k})F_k^{\beta }(1-\phi (\det \sigma _{A}))\big).
\end{eqnarray*}%
We estimate $I$.  Notice that $\phi (\det \sigma _{A})\neq 0$ implies that
$\det \sigma _{A}\geq \varepsilon _{\ast }/2$. We use the integration by parts
with respect to $A$ in Proposition \ref{IBP}, and we obtain
\begin{equation*}
I=\E\big(f_{\delta }(A+B_{k})F_k^{\alpha }H_{\alpha}^{r}(A,\phi (\det
\sigma _{A}))\big).
\end{equation*}%
The estimate (\ref{gamma11}) for the weight gives that
$$
|H^r_\alpha(A,\det\sigma_A)|\leq C\big(1\vee (\det\sigma_A)^{-1}\big)^{r(r+1)}\big(1+|A|_{1,r+1}^{2N(r+2)}+|LA|_{r-1}^2\big)^r\times |\phi(\det\sigma_A)|_r,
$$
$C$ denoting a universal constant.  Since $\sigma_A=\frac 12\sigma_{S_n/2}=\frac 1{n}\sum_{k=1}^{n/2}\chi_k I$, all the Malliavin derivatives are null, so $|\phi(\det\sigma_A)|_r=|\phi(\det\sigma_A)|\leq 1$, so that
$$
|H^r_\alpha(A,\det\sigma_A)|\leq C\big(1\vee (\det\sigma_A)^{-1}\big)^{r(r+1)}\big(1+|A|_{1,r+1}^{2N(r+2)}+|LA|_{r-1}^2\big)^r.
$$
We pass now to expectation: by using the H\"older inequality, we may find some universal constants $C,q,p$ such that
\begin{equation*}
\E\big(| H_{\alpha }^{r}(A,\phi (\det \sigma _{A}))|^{2}\big)
\leq \frac{C}{\varepsilon _{\ast }^{q}}\big(1+\|A\|_{1,r+1,p}+\| LA\|_{r-1,p}\big)^{q}\leq 
C',
\end{equation*}%
the latter inequality following from (\ref{abs18}) and (\ref{abs19}).  Now, $F_k$ and $H_{\alpha}^{r}(A,\phi (\det \sigma _{A}))$
are independent, so that
\begin{equation*}
\left\vert I\right\vert
\leq C\left\Vert f\right\Vert _{\infty }\E(|F|^m).
\end{equation*}%
We estimate now $J.$ By recalling again that $F_k$ and $\sigma _{A}$ are
independent and by using (\ref{abs26}) and (\ref{abs29}), we obtain%
\begin{eqnarray*}
| J|  &\leq &
\|\partial _{\alpha }f_{\delta }\| _{\infty }
\E(|F_k^{\beta}(1-\phi (\det \sigma _{A})| )
\leq \delta ^{-r/2}\|f\|_{\infty }
\E(|F_k|^m)\,\P(\sigma _{S_{n/2}}\leq
\varepsilon _{\ast })) \\
&\leq &C\delta ^{-r/2}\|f\| _{\infty }\E(|F|^m)\times e^{-n/C}.
\end{eqnarray*}%
By resuming, we get
$$
\Big|\E\Big(\partial _{\alpha }f_{\delta }\Big(\widehat{S}%
_{n}^{k}+\lambda \frac{F_k}{n^{1/2}}\Big)F_k^{\beta }\Big)\Big|
\leq C\|f\| _{\infty }\E(|F|^m)\big(1+\delta ^{-r/2} e^{-n/C}\big).
$$
And similarly, we prove that
$$
\Big|\E\Big(\partial _{\alpha }f_{\delta }\Big(\widehat{S}%
_{n}^{k}+\lambda \frac{G_{k}}{n^{1/2}}\Big)G_k^{\beta }\Big)\Big|
\leq C\|f\| _{\infty }\E(|G|^m)\big(1+\delta ^{-r/2} e^{-n/C}\big).
$$
$\square $

\medskip

We can now give a nice estimate for $\mU_r^nf_\delta$ in terms of $\|f\|_\infty$. And this is enough for the moment.

\begin{lemma}\label{4.4bis}
Suppose that $\mu_F\succeq \leb_N$. Let $r\geq 2$ and  $F\in L^{r+1}(\Omega)$. For $f\in L^\infty(\R^N)$ and $\delta>0$, set $f_\delta=f\ast \gamma_\delta$, $\gamma_\delta$ being the centred normal density in $\R^N$ with covariance matrix $\delta I$. Then there exists $C>0$ depending on $r$ and $N$ only such that for every $f\in L^{\infty}(\R^N)$ one has
\begin{equation}\label{Rnrdelta-est}
|\mU^n_r f_\delta|
\leq
C (1+\E(|F|^{r+1}))^{[r/3]\vee 1}
\|f\|_\infty(1+\delta^{-\frac{r+4}2}e^{-n/C})\big(
\sup_{|\alpha|\leq r}|\Delta_\alpha|
\times n^{\frac {r-[r/3]}2}+ n\Big).
\end{equation}
\end{lemma}

\textbf{Proof.}
By using (\ref{TR}) and (\ref{Phikr}), we can write
$$
\mU^n_r f_\delta
=\sum_{k=1}^n\Big[\sum_{j=1}^{k-1}
\sum_{t=3}^r\mathcal{R}^{k-j}_{r-t,n}\Psi^{(j)}_tf_\delta+\mathcal{R}_{r,n}^kf_\delta\Big].
$$
Since $g\mapsto \mathcal{R}^{l}_{t,n}g$ is linear, by using the expansion of $\Psi^{(k)}$ in Lemma \ref{miracle-real} and by recalling that $Q_{i-1}(k)\geq 0$, we get
$$
|\mU^n_r f_\delta|
\leq \sum_{k=2}^n\sum_{j=1}^{k-1}\sum_{t=3}^r
\sum_{i=1}^{[t/3]}Q_{i-1}(j)
|\mathcal{R}^{k-j}_{r-t,n}\mA^i_tf_\delta|+\sum_{k=1}^n|\mathcal{R}_{r,n}^kf_\delta|.
$$
Since $r\geq 2$, (\ref{R01}) gives $\mathcal{R}^{\ell}_{0,n}=\frac 1n \mathcal{R}^{\ell}_{2,n}$ and
$\mathcal{R}^{\ell}_{1,n}=\frac 1{\sqrt n} \mathcal{R}^{\ell}_{2,n}$. So,
we isolate in the sum the terms with $t=r-1,r$ and we obtain
\begin{equation}\label{app}
\begin{array}{ll}
|\mU^n_r f_\delta|
\leq
&\displaystyle\sum_{k=2}^n\sum_{j=1}^{k-1}\Big[\I_{r\geq 5}\sum_{t=3}^{r-2}
\sum_{i=1}^{[t/3]}Q_{i-1}(j)
|\mathcal{R}^{k-j}_{r-t,n}\mA^i_tf_\delta|
+\I_{r\geq 4}\frac 1{\sqrt n}\sum_{i=1}^{[(r-1)/3]}Q_{i-1}(j)
|\mathcal{R}^{k-j}_{2,n}\mA^i_{r-1}f_\delta|\smallskip\\
&\qquad\qquad+\displaystyle\I_{r\geq 3}\frac 1n\sum_{i=1}^{[r/3]}Q_{i-1}(j)
|\mathcal{R}^{k-j}_{2,n}\mA^i_rf_\delta|
\Big]
+\sum_{k=1}^n|\mathcal{R}_{r,n}^kf_\delta|.
\end{array}
\end{equation}
We have (recall formula (\ref{Rkf}))
\begin{align*}
|\mathcal{R}_{r,n}^kf_\delta|
\leq&
n^{\frac{r+1}2}\Big[\Big|\E\Big(U_{r}f_\delta\Big(\widehat S^k_n,\frac{F_k}{\sqrt n}\Big)\Big)\Big|
+\Big|\E\Big(U_{r}f_\delta\Big(\widehat S^k_n,\frac{G_k}{\sqrt n}\Big)\Big)\Big|\Big]\\
&+\sum_{p=3}^r
\!\!\!\!\sum_{\mbox{\scriptsize{$\begin{array}{c}|\alpha|=p\\|\beta|={2[\frac{r-p}2]+2}\end{array}$}}}\!\!\!\!\!\!\!\!
|\Delta_\alpha|
\int_0^1
\Big|\E\Big(\partial_\beta\partial_\alpha f_\delta\Big(\widehat S^{k}_n+\sqrt s\,\frac {G_k}{\sqrt n}\Big)\Big)\Big|ds,
\end{align*}
%
and by using Lemma \ref{4.4-new} we get
\begin{equation}\label{app1}
|\mathcal{R}^k_{r,n} f_\delta|
\leq
C(1+\E(|F|^{r+1}))\|f\|_\infty(1+\delta^{-\frac{r+2}2}e^{-n/C}).
\end{equation}
As for the other sums in the right hand side of (\ref{app}), for $s\geq 2$ we have
$$
|\mathcal{R}^{k-j}_{s,n} \mA^i_tf_\delta|
\leq
\sum_{|\gamma|=t}|c^i_\gamma|\times|\mathcal{R}^{k-j}_{s,n} \partial_\gamma f_\delta|
\leq
C\sup_{|\alpha|\leq t}|\Delta_\alpha|(1+\E(|F|^t))^{i-1}\sum_{|\gamma|=t}|\mathcal{R}^{k-j}_{s,n} \partial_\gamma f_\delta|,
$$
last inequality following from (\ref{ci-est}). We use again Lemma \ref{4.4-new}: for $|\gamma|=t$,
$$
|\mathcal{R}^{k-j}_{s,n} \partial_\gamma f_\delta|
\leq
C(1+\E(|F|^{s+1}))\|f\|_\infty(1+\delta^{-\frac{s+t+2}2}e^{-n/C}).
$$
We apply such inequality with: $t\leq r-2$ and $s=r-t$, $t=r-1$ and $s=2$, $t=r$ and $s=2$. Then,
\begin{align*}
|\mathcal{R}^{k-j}_{r-t,n} \mA^i_tf_\delta|
\leq
&C\sup_{|\alpha|\leq r}|\Delta_\alpha|(1+\E(|F|^{r+1}))^{[r/3]\vee 1}\|f\|_\infty(1+\delta^{-\frac{r+2}2}e^{-n/C})\\
|\mathcal{R}^{k-j}_{2,n} \mA^i_{r-1}f_\delta|
\leq
&C\sup_{|\alpha|\leq r}|\Delta_\alpha|(1+\E(|F|^{r+1}))^{[r/3]\vee 1}\|f\|_\infty(1+\delta^{-\frac{r+3}2}e^{-n/C})\\
|\mathcal{R}^{k-j}_{2,n} \mA^i_rf_\delta|
\leq
&C\sup_{|\alpha|\leq r}|\Delta_\alpha|(1+\E(|F|^{r+1}))^{[r/3]\vee 1}
\|f\|_\infty(1+\delta^{-\frac{r+4}2}e^{-n/C})
\end{align*}
By inserting such estimates and (\ref{app1}) in (\ref{app}), we get
\begin{align*}
|\mU^n_r f_\delta|
\leq
&\sum_{k=2}^n\sum_{j=1}^{k-1}\Big[\I_{r\geq 5}\sum_{t=3}^{r-2}
\sum_{i=1}^{[t/3]}Q_{i-1}(j)
|\mathcal{R}^{k-j}_{r-t,n}\mA^i_tf_\delta|
+\I_{r\geq 4}\frac 1{\sqrt n}\sum_{i=1}^{[(r-1)/3]}Q_{i-1}(j)
|\mathcal{R}^{k-j}_{2,n}\mA^i_{r-1}f_\delta|\smallskip\\
&\qquad\qquad+\I_{r\geq 3}\frac 1n\sum_{i=1}^{[r/3]}Q_{i-1}(j)
|\mathcal{R}^{k-j}_{2,n}\mA^i_rf_\delta|
\Big]
+\sum_{k=1}^n|\mathcal{R}_{r,n}^kf_\delta|\\
\leq
&C (1+\E(|F|^{r+1}))^{[r/3]\vee 1}
\|f\|_\infty(1+\delta^{-\frac{r+4}2}e^{-n/C})\times\\
&\times\Big(
\sup_{|\alpha|\leq r}|\Delta_\alpha|
\sum_{k=2}^n\sum_{j=1}^{k-1}\Big[\I_{r\geq 5}\!\!\!
\sum_{i=1}^{[(r-2)/3]}Q_{i-1}(j)+\I_{r\geq 4}\frac 1{\sqrt n}\!\!\!\sum_{i=1}^{[(r-1)/3]}Q_{i-1}(j)+\\
&\qquad
+\I_{r\geq 3}\frac 1n\!\sum_{i=1}^{[r/3]}Q_{i-1}(j)
\Big]+n\Big)
\end{align*}
Since $\sum_{k=2}^n\sum_{j=1}^{k-1}
\sum_{i=1}^{L}Q_{i-1}(j)=\sum_{i=1}^{L}P_{i+1}(n-1)$ is a polynomial of order $L+1$ we obtain
\begin{align*}
|\mU^n_r f_\delta|
\leq
&C (1+\E(|F|^{r+1}))^{[r/3]\vee 1}
\|f\|_\infty(1+\delta^{-\frac{r+4}2}e^{-n/C})\times\\
&\times\Big(
\sup_{|\alpha|\leq r}|\Delta_\alpha|
\Big[\I_{r\geq 5}
n^{[(r-2)/3]+1}+\I_{r\geq 4} n^{[(r-1)/3]+\frac 12}
+\I_{r\geq 3} n^{[r/3]}
\Big]+n\Big)
\end{align*}
and the statement follows by noticing that
$$n^{[(r-2)/3]+1}\I_{r\geq 5}+n^{[(r-1)/3]+1/2}\I_{r\geq 4}+n^{[r/3]}\I_{r\geq 3}
\leq C n^{[r/3]+\frac {r-3[r/3]}2}.
$$
$\square$

\subsection{Estimate of the error in total variation distance}\label{sect-errTV}

We want to get rid of the derivatives of $f$ which appear in the coefficients $\mD_m f$. In order to do it we will use integration by
parts w.r.t. the Gaussian law and then the Hermite polynomials come on. Again, we assume $\mu_F\succeq \leb_N$ and  $F$ has null mean and identical covariance matrix.

We denote by $H_{m}$ the Hermite polynomial of order $m$ on $\R$, that is,
\begin{equation}
H_{m}(x)=(-1)^{m}e^{\frac{1}{2}x^{2}}\frac{d^{m}}{x^{m}}e^{-\frac{1}{2}x^{2}}.
\label{tv1}
\end{equation}%
For a multiindex $\alpha=(\alpha _{1},...,\alpha _{r})\in \{1,...,N\}^{r}$ we
denote $\beta _{i}(\alpha )=\card\{j:\alpha _{j}=i\}$ so that $\partial
_{\alpha }=\partial _{x_{1}}^{\beta _{1}(\alpha )}\ldots\partial
_{x^{N}}^{\beta _{d}(\alpha )}.$ And we define the Hermite polynomial on $\R^{N}
$ corresponding to the multiindex $\alpha $ by%
\begin{equation}
H_{\alpha }(x)=\prod_{i=1}^{N}H_{\beta _{i}(\alpha )}(x_{i})\qquad \mbox{for}\qquad
x=(x_{1},...,x_{N}).  \label{tv2}
\end{equation}%
With this definition we have
\begin{equation*}
\partial _{\alpha }e^{-\frac{1}{2}\left\vert x\right\vert
^{2}}=(-1)^{\left\vert \alpha \right\vert }H_{\alpha }(x)e^{-\frac{1}{2}%
\left\vert x\right\vert ^{2}}
\end{equation*}%
and using integration by parts, for a centred Gaussian random variable $G\in
\R^{N}$%
\begin{equation}
\E(\partial _{\alpha }f(G))=\E(f(G)H_{\alpha }(G)).  \label{tv3}
\end{equation}%
This means that we can compute $\E(\mA^i_tf(G))$ by means of $f$ and not of its derivatives. In fact, for $i\geq 1$ and $t\geq 0$, we define the polynomials $\mH^i_t(x)$ as follows:
\begin{equation}\label{Hi}
\mH^{i}_t(x)=\sum_{\alpha : |\alpha|=t}c^i_\beta H_{\alpha}(x),\quad\mbox{$c^i_\alpha$ defined in (\ref{ci}) and $H_\alpha$ given in (\ref{tv2})}.
\end{equation}
Since $\mA^i_t=\sum_{\alpha\,:\,|\alpha|=t}c^i_\alpha\partial_\alpha$, (\ref{tv3}) gives
$$
\E(\mA^i_tf(G))=\sum_{\alpha\,:\,|\alpha|=t}c^i_\gamma\E(\partial_\alpha f(G))
=\sum_{\alpha\,:\,|\alpha|=t}c^i_\alpha \E(f(G)H_\alpha(G))=\E(f(G)\mH^i_t(G)).
$$
Therefore, for every $f\in C^m_b(\R^N)$ the coefficients $\mD_m f$, $m\geq 1$, in (\ref{D}) can be written as
\begin{equation}\label{Dbis}
\begin{array}{l}
\displaystyle
\mD_m f=\E(f(G)\mK_m(G)),\quad m\geq 1,\quad\mbox{where}\smallskip\\
\displaystyle
\mK_m(x)
=
\sum_{\mbox{\scriptsize{$\begin{array}{c}t=3\vee m\\t-m \mbox{ even}\end{array}$}}}^{3m}
\sum_{i=1\vee \frac{t-m}2}^{[t/3]}a_{i,\frac{t-m}2}
\mH^i_t(x),\quad \mbox{$a_{i,l}$ given in (\ref{aip}) and $\mH^i_t$ given in (\ref{Hi})}.
\end{array}
\end{equation}

We are now ready to tackle our original problem: the exact expansion in total variation distance of the law $\mu_n$ of $S_n$. To this purpose, for $r\geq 2$ and $n\geq 1$  we define the following measure in $\R^N$:
\begin{equation}\label{Gammanp}
\Gamma_{n,r}(dx)=\gamma(x)\Big(1+\sum_{m=1}^{[r/3]} \frac 1{n^{\frac m2}}\mK_m(x)\Big)dx,\quad\mbox{$\mK_m(x)$ given in (\ref{Dbis})},
\end{equation}
where $\gamma(x)$ denotes the probability density function of a standard normal random variable in $\R^N$.
We stress that $\Gamma_{n,r}(dx)=\gamma(x)dx=:\Gamma(dx)$ not only for $r= 2$ but also when $\Delta_\alpha=0$ for every $|\alpha|\leq r$. In fact, in  the latter case, (\ref{ci-est}) gives $c^i_\alpha=0$ for every $i\geq 1$ and $|\alpha|\leq r$,
then from  (\ref{Hi}) we have $\mathcal{H}^i_t\equiv 0$ for every $i\geq 1$ and $t\leq r$ and from (\ref{Dbis}) we obtain $\mathcal{K}_m\equiv 0$ for every $m\leq r$.

\begin{theorem}\label{main-th}
Suppose $\mu_F\succeq \leb_N$. Let $r\geq 2$ and $F\in L^{r+1}(\Omega)$. For $n\geq 1$, let $\mu_n$ denote the law of $S_n$  and $\Gamma_{n,r}$ stand for the measure in (\ref{Gammanp}).
Then
there exists a constant $C>0$ depending on $r$ and $N$ only such that for every $n\in\N$,
$$
d_{TV}(\mu_n,\Gamma_{n,r})\leq C \big(1+\E(|F|^{r+1})\big)^{[r/3]\vee 1}\Big[\sup_{|\alpha|\leq r}|\Delta_\alpha|\times \frac{1}{n^{\frac{[r/3]+1}2}}+\frac 1{n^{\frac{r-1}2}}\Big].
$$
\end{theorem}

\textbf{Proof}. We study $|\int fd\mu_n-\int fd\Gamma_{n,r}|$ for $f\in L^{\infty}(\R^N)$. From now on, $C$ will denote a constant, possibly varying from line to line, that may depend only on $N$ and $r$.

\smallskip

We take $\delta >0$ and we consider the regularized function $f_{\delta }=f\ast
\gamma _{\delta }$ where $\gamma _{\delta }$ is the centred Gaussian density of
covariance matrix $\delta I$. We have
$$
\Big|\int fd\mu_n-\int fd\Gamma_{n,r}\Big|
\leq I_{n,\delta}+I'_{n,\delta}+J_{n,\delta}
$$
with
$$
I_{n,\delta}=\Big|\int (f-f_\delta)d\mu_n\Big|,\quad
I'_{n,\delta}=\Big|\int (f-f_\delta)d\Gamma_{n,r}\Big|,\quad
J_{n,\delta}=
\Big|\int f_\delta d\mu_n-\int f_\delta d\Gamma_{n,r}\Big|
$$
By (\ref{abs30'})
$$
I_{n,\delta}
\leq
C\left\Vert f\right\Vert _{\infty }(1+\E(|F|))\big(e^{-n/C}+\delta ^{1/b} n^{(b-2)/(2b)}\big),
$$
where $b>4$ is a suitable constant, independent of $F$ and $f$. And using standard integration by parts on $\R^{N}$,
$$
I'_{n,\delta}
\leq
C\left\Vert f\right\Vert _{\infty }\delta ^{1/2}.
$$
Moreover, since
$$
\int f_\delta d\Gamma_{n,r}
=\E(f_\delta(G))+\sum_{m=1}^{[r/3]}\frac 1{n^{\frac m2}}\mD_m f_\delta,
$$
Theorem \ref{dev1} gives
$$
J_{n,\delta}
=\frac 1{n^{\frac{[r/3]+1}2}}|\mE^n_{r} f_\delta|
$$
with
\begin{align*}
|\mE^n_{r} f_\delta|
\leq
&n^{\frac{[r/3]+1}2} \Big[
\sum_{m=[r/3]+1}^r\frac 1{n^\frac m2}\sum_{t=3\vee m}^{(3m)\wedge r}\sum_{i=1\vee\frac {t-m}2}^{[t/3]}|a_{i,\frac {t-m}2}|\times|\E(\mA^i_tf_\delta(G))|+\frac 1{n^{\frac{r+1}2}}|\mU^n_r f_\delta|\Big]
\end{align*}
But since $\E(\mA^i_tf_\delta(G))=\E(f_\delta(G)\mH^i_t(G))$, then
$$
|\E(\mA^i_tf_\delta(G))|\leq \|f_\delta\|_{\infty}\E(|\mH^i_t(G)|)\leq C \|f\|_\infty(1+\E(|F|^{t-1}))\sup_{|\alpha|\leq t}|\Delta_\alpha|.
$$
We use now Lemma \ref{4.4bis}: for $r\geq 2$, we apply (\ref{Rnrdelta-est}) and we get
\begin{align*}
|\mE^n_{r} f_\delta|
&\leq
C(1+\E(|F|^{r+1}))^{[r/3]\vee 1}\|f\|_{\infty}(1+\delta^{-\frac{r+4}2}e^{-n/C})
\Big[\sup_{|\alpha|\leq r}|\Delta_\alpha|+\frac 1{n^{\frac{r-[r/3]-2}2}}\Big].
\end{align*}
By replacing, we get
$$
J_{n,\delta}\leq
C\,(1+\E(|F|^{r+1}))^{[r/3]\vee 1}\|f\|_{\infty}(1+\delta^{-\frac{r+4}2}e^{-n/C})
\Big[\sup_{|\alpha|\leq r}|\Delta_\alpha|\times \frac 1{n^{\frac{[r/3]+1}2}}+\frac 1{n^{\frac{r-1}2}}\Big].
$$
By resuming, we can write
\begin{align*}
\Big|\int fd\mu_n-\int fd\Gamma_{n,r}\Big|
\leq
&
C\left\Vert f\right\Vert _{\infty }(1+\E(|F|^{r+1}))^{[r/3]\vee 1}
\Big[e^{-\frac{n}{C}}+\delta ^{1/2}+\delta ^{1/b} n^{(b-2)/(2b)}+\\
&\qquad+(1+\delta^{-\frac{r+4}2}e^{-n/C})
\Big(\sup_{|\alpha|\leq r}|\Delta_\alpha|\times \frac 1{n^{\frac{[r/3]+1}2}}+\frac 1{n^{\frac{r-1}2}}\Big)\Big].
\end{align*}
Now, we choose $\delta=\delta_n$ such that $\delta_n ^{1/b} n^{(b-2)/(2b)}=\frac 1{n^{\frac{r-1}2}}$. By observing that $n\mapsto \delta_n^{-\frac{r+4}2}e^{-n/C}$ is bounded and $\delta_n^{1/2}\leq \frac 1{n^{\frac{r-1}2}}$, we get
$$
\Big|\int fd\mu_n-\int fd\Gamma_{n,r}\Big|
\leq
C\left\Vert f\right\Vert _{\infty }(1+\E(|F|^{r+1}))^{[r/3]\vee 1}
\Big[\sup_{|\alpha|\leq r}|\Delta_\alpha|\times\frac 1{n^{\frac{[r/3]+1}2}}+\frac 1{n^{\frac{r-1}2}}\Big]
$$
and the result follows.
$\square$

\medskip

We can now pass to the following.

\medskip

\textbf{Proof of Theorem \ref{Speed-r} and \ref{main-th-1}.} We apply Theorem \ref{main-th} with $F$ replaced by $A(F)F$, where $A(F)$ is the inverse of $C(F)^{1/2}$, $C(F)$ denoting the covariance matrix. And it is clear that now the constants appearing in the estimates  will depend on  $C(F)$ as well, through its most significant eigenvalues (the smallest and the largest one; see, e.g.,  Remark \ref{noId1} and \ref{noId2}). $\square$

\medskip

We conclude by explicitly writing $\mK_m(x)$ for $m=1,2,3$.
From (\ref{Dbis}) we have:
\begin{align*}
&\mK_1(x)=a_{1,1}\mH^1_3(x)\\
&\mK_2(x)=a_{1,1}\mH^1_4(x)+a_{2,2}\mH^2_6(x)\\
&\mK_3(x)=a_{1,0}\mH^1_3(x)+a_{1,1}\mH^1_5(x)+a_{2,2}\mH^2_7(x)+a_{3,3}\mH^3_9(x),
\end{align*}
where $\mH^i_t(x)=\sum_{|\gamma|=t}c^i_\gamma H_\gamma(x)$. Now, from (\ref{ci}) it is easy to see that
\begin{align*}
&c^1_\gamma=\left\{
\begin{array}{ll}
\displaystyle
\frac 1{3!}\Delta_\gamma & \mbox{if $|\gamma|=3$}\smallskip\\
\displaystyle
\frac 1{4!}\Delta_\gamma & \mbox{if $|\gamma|=4$}\smallskip\\
\displaystyle
\frac 1{5!}\Delta_\gamma
-\frac1 {3!2!}\Delta_{(\gamma_1,\gamma_2,\gamma_3)}1_{\gamma_4=\gamma_5}& \mbox{if $|\gamma|=5$},
\end{array}
\right.\\
&c^2_\gamma=\left\{
\begin{array}{ll}
\displaystyle
\frac 1{(3!)^2}\Delta_{(\gamma_1,\gamma_2,\gamma_3)}\Delta_{(\gamma_4,\gamma_5,\gamma_6)}& \mbox{if $|\gamma|=6$}\smallskip\\
\displaystyle
\frac 1{3!4!}\big(\Delta_{(\gamma_1,\gamma_2,\gamma_3)}\Delta_{(\gamma_4,\gamma_5\gamma_6,\gamma_7)}
+\Delta_{(\gamma_1,\gamma_2,\gamma_3,\gamma_4)}\Delta_{(\gamma_5,\gamma_6,\gamma_7)}\big)& \mbox{if $|\gamma|=7$},
\end{array}
\right.\\
&c^3_\gamma=
\frac 1{(3!)^3}\Delta_{(\gamma_1,\gamma_2,\gamma_3)}
\Delta_{(\gamma_4,\gamma_5,\gamma_6)}
\Delta_{(\gamma_7,\gamma_8,\gamma_9)}
\quad \mbox{if $|\gamma|=9$}.
\end{align*}
Moreover, $a_{1,0}=0$, $a_{1,1}=1$, $a_{2,2}=b_{1,2}=\frac 12 B_0=\frac 12$ and $a_{3,3}=a_{2,2}b_{2,3}=\frac 12\cdot\frac 13 B_0=\frac 16$. So, we can write
\begin{align*}
\mK_1(x)=&\frac {1}{3!}\sum_{|\gamma|=3}\Delta_\gamma H_\gamma(x)\\
\mK_2(x)=&\frac {1}{4!}\sum_{|\gamma|=4}\Delta_\gamma H_\gamma(x)+
\frac {1}{2(3!)^2}\sum_{|\gamma|=6}\Delta_{(\gamma_1,\gamma_2,\gamma_3)}\Delta_{(\gamma_4,\gamma_5,\gamma_6)}H_\gamma(x)\\
\mK_3(x)=&\sum_{|\gamma|=5}\Big(
\frac 1{5!}\Delta_\gamma-\frac1 {2\times3!}\Delta_{(\gamma_1,\gamma_2,\gamma_3)}1_{\gamma_4=\gamma_5}\Big)\,H_\gamma(x)+\\
&+\frac 1{2\times 3!4!}\sum_{|\gamma|=7}
\big(\Delta_{(\gamma_1,\gamma_2,\gamma_3)}\Delta_{(\gamma_4,\gamma_5,\gamma_6,\gamma_7)}
+\Delta_{(\gamma_1,\gamma_2,\gamma_3,\gamma_4)}\Delta_{(\gamma_5,\gamma_6,\gamma_7)}\big)
H_\gamma(x)+\\
&+\frac {1}{6\times (3!)^3}\sum_{|\gamma|=9}\Delta_{(\gamma_1,\gamma_2,\gamma_3)}
\Delta_{(\gamma_4,\gamma_5,\gamma_6)}
\Delta_{(\gamma_7,\gamma_8,\gamma_9)}
H_\gamma(x)
\end{align*}
In the case $N=1$, for $t\in\N$ set
$$
\ell_t=\frac{\E(F^t)}{\Var(F)^{t/2}}.
$$
Note that $\ell_t$ is strictly connected to the Lyapunov ratio  $L_t=\frac{\E(|F|^t)}{\Var(F)^{t/2}}$. By recalling that
for $G\sim\mathrm{N}(0,1)$ then $\E(G^t)=0$ if $t$ is odd and $\E(G^t)=(t-1)!!$ if $t$ is even (with the convention $(-1)!!=1$), we obtain
$\Delta_t=\ell_t$ if $t$ is odd and $\Delta_t=\ell_t-(t-1)!!$ if $t$ is even. Remark that $\Delta_3=\ell_3$ and $\Delta_4=\ell_4-3$ are the skewness and  the kurtosis respectively. Hence, we obtain the polynomials in the classical  Edgeworth expansion:
\begin{align*}
\mK_1(x)=\frac {\ell_3}{6}H_3(x),\qquad
\mK_2(x)=\frac {(\ell_4-3)}{24} H_4(x)+
\frac {\ell_3^2}{72}H_6(x)\\
\mK_3(x)=
\Big(\frac{\ell_5}{5!}-\frac{\ell_3}{2\times 3!}\Big)H_5(x)
+\frac{\ell_3(\ell_4-3)}{3!4!}H_7(x)+\frac{\ell_3^3}{6(3!)^3}H_9(x).
\end{align*}

\appendix
\section{Probability measures which are locally lower bounded by the Lebesgue me\-asure}\label{app-1}

We discuss here the proof of Proposition \ref{prop-main1}.
For a random variable $F\in \R^{N}$ with law $\mu_F$, we recall that $\mu_F\succeq\leb_N$ if there exists an open set $D\subset \R^{N}$ and $\varepsilon >0$ such that%
\begin{equation}
\mu _{F}(A):=\P(F\in A)\geq \varepsilon \leb _{N}(A\cap D)\qquad \forall
A\in \mathcal{B}(\R^{N}).  \label{AP1}
\end{equation}%

Remark that we have already proved that if $\mu_F\succeq \leb_N$ then (\ref{abs7}) holds (see Proposition \ref{prop-chi}).

We first prove the equivalence $(i)\Leftrightarrow (ii)$:
\begin{lemma}\label{lemma-app-1}
$\mu _{F}\succeq \leb _{N}$ if and only if there exists a non-negative
measure $\mu $ with $\mu (\R^{N})<1$ and a non-negative lower semicontinuous
function $p$ with $\int_{\R^{N}}p(v)dv=1-\mu (\R^{N})$ such that
\begin{equation}
\mu_F(dv)=\mu (dv)+p(v)dv.  \label{Condition}
\end{equation}
\end{lemma}

\textbf{Proof}. If (\ref{AP1}) holds we take $v_0\in D$ and $r>0$ such that $B_r(v_0)\subset D$. Then, it suffices to take $p(x)=\varepsilon
1_{B_{r}(v_{0})}(x)$ and $\mu (A)=\P(F\in A)-\int_{A}p(v)dv$, which turns out to be a non-negative measure.

Suppose now that (\ref{Condition}) holds. Since $p$ is non-negative and lower semicontinuous we may find an increasing sequence of non-negative and continuous functions $%
p_{n},n\in \N$ such that $p_{n}\uparrow p.$ It follows that $\int
p_{n}\uparrow \int p=1-\mu (\R^{N})>0$, and  we may find $n$ such that $\int
p_{n}>0$. So there exists $v_{0}$ such that $p_{n}(v_{0})>0.$ Since $p_{n}$
is continuous, this implies that $p(v)\geq p_{n}(v)\geq \frac{1}{2}%
p_{n}(v_{0})$ for $\left\vert v-v_{0}\right\vert <r$ for some small $r.$ $%
\square $

\medskip

As a consequence we get
the final property in Proposition \ref{prop-main1}:

\begin{lemma}\label{lemma-app-2}
If $\mu _{F}\succeq \leb _{N}$, then the covariance matrix of $F$ is
invertible.
\end{lemma}

\textbf{Proof}. We fix $v_{0}\in \R^{N}$ and $\varepsilon >0$ such that (\ref%
{AP1})\ holds with $D=B_{r}(v_{0}).$ We assume that $\E(F^{i})=0$ so that the
covariance matrix is given by $C^{i,j}(F)=\E(F^{i}F^{j}).$ Then, for $\xi \in
\R^{N}$ we write%
\begin{equation*}
\left\langle C(F)\xi ,\xi \right\rangle =\E(\left\langle F,\xi \right\rangle
^{2})\geq \varepsilon \int_{B_{r}(v_{0})}\left\langle v,\xi \right\rangle
^{2}dv.
\end{equation*}%
We denote $A_{\delta }(\xi )=\{v:\left\langle v,\xi \right\rangle ^{2}\geq
\delta \left\vert \xi \right\vert ^{2}\}$ and we note that we may choose $%
\delta (v_{0},r)$ such that
\begin{equation*}
\inf_{\left\vert \xi \right\vert =1}\leb_N (A_{\delta (v_{0},r)}(\xi
))=:\eta (v_{0},r)>0.
\end{equation*}%
Then
\begin{equation*}
\inf_{\left\vert \xi \right\vert =1}\left\langle C(F)\xi ,\xi \right\rangle
\geq \varepsilon \eta (v_{0},r)\leb_N (B_{r}(v_{0})).
\end{equation*}
$\square$

\medskip

We have already proved in Proposition \ref{prop-chi} the implication $(i)\Rightarrow (iii)$. Last implication $(iii)\Rightarrow (ii)$ is trivial. In fact, let
$$
\P(\chi V+(1-\chi )W\in dv)=\P(F\in dv)
$$
where $\chi$ is a Bernoulli r.v. with parameter $p>0$, $V$ in $\R^N$ is absolutely continuous and $W$ is a r.v. in $\R^{N}$. Setting $\mu_F$, $\mu_V$ and $\mu_W$, the law of $F$, $V$ and $W$, respectively, then
$$
\mu_F(dv)=p\mu_V(v)dv+(1-p)\mu_W(dv),
$$
so $F$ has an absolutely continuous component.

\section{Estimates for the Sobolev norms in Lemma \ref{pippo}}\label{app-sob}

This section is devoted to the proof the estimates used in Lemma \ref{pippo}, that is the following.

\begin{lemma}\label{estK}
Let $d\geq 1$, $m\in\N$, $p\geq 1$. Then there exists $C>0$ such that for every $K>1$ and $X=(X^1,\ldots,X^d)$ the following estimates holds:
\begin{align}
&\|\Psi_K(X)X\|_{m,p}
\leq
CK\big(1+\|X\|_{1,m,(m+1)p}\big)^{m+1},\label{estK-1}\\
&
\|L(\Psi_K(X)X)\|_{m,p}
\leq CK\big(1+\|X\|_{1,m+1,4(m\vee 2)p}\big)^{2m+3}\big(1+
\|LX\|_{m,4p}\big)\label{estK-2}
\end{align}
where $\Psi_K(X)$ denote any function in $C^\infty(\R^d)$ such that $1_{B_{K}(0)}\leq \Psi_K\leq 1_{B_{K+1}(0)}$ and whose derivatives are uniformly bounded, that is there exists $L>0$ such that $|\partial_\alpha\Psi_K|\leq L$ for every multiindex $\alpha$.
\end{lemma}

\textbf{Proof.}
For a multiindex $\alpha$, one has
$$
D_\alpha(\Psi_K(X)X^i)
=D_\alpha\Psi_K(X)X^i+\sum_{\beta,\gamma\in A_\alpha, |\beta|\geq 1}
D_\gamma\Psi_K(X)D_\beta X^i
$$
where the condition ``$\beta,\gamma\in A_\alpha$'' means that
$\beta,\gamma$ is a partition of $\alpha$. Moreover, one has
$$
D_\gamma \Psi_K(X)
=\sum_{\ell=1}^{|\gamma|}\sum_{|\rho|=\ell}\partial_\rho\Psi_K(X)\sum_{\beta_1,\ldots,\beta_\ell\in\mathcal{B}_\gamma}
D_{\beta_1}X^{\rho_1}\cdots D_{\beta_\ell}X^{\rho_\ell}
$$
where ``$\beta_1,\ldots,\beta_\ell\in\mathcal{B}_\gamma$'' means that $\beta_1,\ldots,\beta_\ell$ are non-empty multiindexes of $\gamma$ running through the list of all of the (non-empty) ``blocks'' of $\gamma$.
Then, for $|\gamma|\leq m$ we obtain
\begin{equation}\label{est-0}
|D_\gamma \Psi_K(X)|
\leq C\,1_{|X|\leq K+1}\Big(1+\sum_{1\leq|\rho|\leq m} |D_\rho X|\Big)^{m}
\end{equation}
So, for $|\alpha|=m$ we have
\begin{align*}
|D_\alpha(\Psi_K(X)X)|
&\leq
CK\big(1+|X|_{1,m}\big)^{m+1}
\end{align*}
and (\ref{estK-1}) follows.
Consider now $L(\Psi_K(X)X^l)$. We have
\begin{align*}
-L(\Psi_K(X)X^l)
&=-L\Psi_K(X)X^l-\Psi_K(X)LX^l
+\sum_{k=1}^{n}\sum_{i=1}^{d}D_{(k,i)}\Psi_K(X)D_{(k,i)}X^l.
\end{align*}
We use now the inequality $\|XY\|_{m,p}\leq C\|X\|_{m,2p}\|Y\|_{m,2p}$. But concerning the first term of right hand side of the equality above, we take care of the derivatives of $\Psi_K$ as done to obtain formula (\ref{est-0}) and we get
\begin{align*}
\|L(\Psi_K(X)X)\|_{m,p}
&\leq C\|L\Psi_K(X)\|_{m,2p}(\|X1_{|X|<K+1}\|_{2p}+\|X\|_{1,m,2p})\\
&\leq CK \|L\Psi_K(X)\|_{m,2p}(1+\|X\|_{1,m,2p}).
\end{align*}
So, we obtain
\begin{align*}
\|L(\Psi_K(X)X)\|_{m,p}
\leq &C\big( K \|L\Psi_K(X)\|_{m,2p}(1+\|X\|_{1,m,2p})+\\
&+\|\Psi_K(X)\|_{m,2p}\|LX\|_{m,2p}
+\|\Psi_K(X)\|_{1,m,2p}\|X\|_{1,m,2p}\big).
\end{align*}
(\ref{est-0}) gives that
\begin{equation}\label{estK-3}
\|\Psi_K(X)\|_{m,2p}\leq C(1+\|X\|_{1,m,2mp})^m,
\end{equation}
so we can write
\begin{align*}
\|L(\Psi_K(X)X)\|_{m,p}
\leq &CK\big(1+\|X\|_{1,m,2mp}\big)^{m+1}\big(1+ \|L\Psi_K(X)\|_{m,2p}+\|LX\|_{m,2p}\big)
\end{align*}
It remains to estimate $\|L\Psi_K(X)\|_{m,2p}$. Since
$$
L\Psi_K(X)=\sum_{j=1}^d\partial_j\Psi_K(X)LX^j-\frac 12\sum_{i,j=1}^d\partial_i\partial_j\Psi_K(X)\<DX^i,DX^j\>
$$
we have
$$
\|L\Psi_K(X)\|_{m,2p}\leq C\big(\|\nabla\Psi_K(X)\|_{m,4p}\|LX\|_{m,4p}+\|\nabla^2\Psi_K(X)\|_{m,4p}\|DX\|_{m,8p}^2\big).
$$
An inequality analogous to (\ref{estK-3}) can be proved for $\nabla\Psi_K$ and $\nabla^2\Psi_K$, so we obtain
\begin{align*}
\|L\Psi_K(X)\|_{m,2p}
&\leq C\big(
(1+\|X\|_{1,m,4mp})^m\|LX\|_{m,4p}+(1+\|X\|_{1,m,4mp})^m\|X\|_{1,m+1,8p}^2\big)\\
&\leq C(1+\|X\|_{1,m+1,4(m\vee 2)p})^{m+2}\big(1+
\|LX\|_{m,4p}\big).
\end{align*}
Therefore, we can write
\begin{align*}
\|L(\Psi_K(X)X)\|_{m,p}
\leq &CK\big(1+\|X\|_{1,m,2mp}\big)^{m+1}(1+\|X\|_{1,m+1,4(m\vee 2)p})^{m+2}\times\\
 &\quad \times \big( 1+\|LX\|_{m,4p}+\|LX\|_{m,2p}\big)\\
\leq &CK\big(1+\|X\|_{1,m+1,4(m\vee 2)p}\big)^{2m+3}\big(1+
\|LX\|_{m,4p}\big)
\end{align*}
and the statement holds.
$\square$

\section{A backward Taylor formula for the Gaussian law}\label{BT}

We give here a simple result on a Taylor formula of a backward type for the normal law.

\begin{lemma}\label{lemma-Taylor-back}
Let $G$ denote a centred normal distributed r.v. in $\R^N$.
Then for every $L\in\N$ and $g\in C_b^{2(L+1)}(\R^N)$ one has
\begin{equation}\label{Taylor-back}
g(0)=\sum_{\ell=0}^L \frac{(-1)^\ell}{2^\ell \ell!}\sum_{|\beta|={2\ell}}\theta_\beta\E(\partial_\beta g(W_1))
+ \frac{(-1)^{L+1}}{2^{L+1} L!}\sum_{|\beta|={2L+2}}\theta_\beta\int_0^1 s^L
\E(\partial_\beta g(W_s))ds.
\end{equation}
$\theta_\beta$ being defined in (\ref{DO1}).
\end{lemma}
\textbf{Proof.}
Let $W$ denote a Brownian motion in $\R^N$. By It\^o's formula, one has $\E(g(W_1))=g(W_t)+\frac 12\int_t^1\E(\Delta g(W_s))ds$, so we can write
\begin{equation}\label{Tay-app}
\E(g(W_t))=g(W_1)-\frac 12\sum_{|\beta|=2}\theta_\beta\int_t^1\E(\partial_\beta g(W_s))ds.
\end{equation}
Taking $t=0$, this gives $g(0)=\E(g(W_1))-\frac 12\sum_{|\beta|=2}\theta_\beta\int_0^1\E(\partial_\beta g(W_s))ds$ and, by iteration, we write
$$
g(0)=\E(g(W_1))-\frac 12\sum_{|\beta|=2}\theta_\beta\E(\partial_\beta g(W_1))
-\frac 12\sum_{|\beta|=2}\theta_\beta\int_0^1\big[\E(\partial_\beta g(W_s))-\E(\partial_\beta g(W_1))\big]ds.
$$
By using (\ref{Tay-app}) we get
\begin{align*}
g(0)
&=\E(g(W_1))-\frac 12\sum_{|\beta|=2}\theta_\beta\E(\partial_\beta g(W_1))
+\frac 14\sum_{|\beta|=4}\theta_\beta\int_0^1u
\E(\partial_\beta g(W_u))du.
\end{align*}
(\ref{Taylor-back}) now follows by further iteration and by recalling that $W_s$ and $\sqrt sG$ have the same law.
$\square$

\end{document}